
\tolerance=10000
\magnification=\magstep1
\def\nagy{\font\caps=cmcsc10\caps}
\def\kis{\font\kisf=cmr5\kisf}

\def\sethoff{\hoffset=0.5 true in}
\sethoff  
\hsize=6.0 true in
\voffset=.3 true in
\vsize=8.5 true in  

\newif\ifduplexpr
\duplexprfalse

\def\duplex{\hoffset=0.25 true in \duplexprtrue \marginfalse}


\def\uj {\bigskip \rm}


\def\redef#1#2{\expandafter\ifx\csname #1\endcsname\relax
    \expandafter\edef\csname #1\endcsname{#2} 
    \else \message{redefinition of '\string#1'
    } \fi }

\newif\ifexist

\def\testfile#1{
\openin 0=#1
\ifeof 0
\existfalse
\else
\existtrue
\fi
\closein 0
}


\testfile{cite.inc}
\ifexist
\input cite.inc 
\else
\message{!!!!!!!!!! One more pass needed for references !!!!!!!!!}
\fi


\immediate\openout 0=cite.inc


\immediate\openout 2=conten.inc
\def\writeitem#1#2{\write2
  {\string
  \line{#1\string\nagy{}#2\string\rm\string\leaderfill\string\quad\folio}}}


\def\boxit#1{\vbox{\hrule\hbox{\vrule#1\vrule}\hrule}}


\newif\ifmargin
\margintrue


\newcount\sorszam \sorszam=0
\def\sorszaminc {\advance\sorszam by1 
\ifnum \section=0 \else \the\section.\fi
\the\sorszam. }

\def\authorstr{}
\def\shorttitlestr{}

\def\author#1{\def\authorstr{#1} 
 \bigskip
 \centerline{#1}
}

\def\abstract#1{\bigskip
 {\advance\leftskip by 1in
 \centerline{\nagy Abstract}\medskip\noindent
 #1}
}

\def\shorttitle#1{\def\shorttitlestr{#1}}

\def\oldalszam{
\headline={
 \nagy
 \ifnum \folio > 1
 \ifodd \folio {\hfil \shorttitlestr \hfil\folio} 
 \else
 {\folio\hfil \authorstr \hfil}
 \fi
 \else \hfil
 \fi
}
}

\nopagenumbers
\oldalszam


\newcount\section \section=0
\def\newsection#1{
      \goodbreak
      \advance\section by1 \sorszam=0 
      \medskip\bigskip\centerline{\nagy \the\section. #1}
      \nobreak\nobreak
     \writeitem{\the\section.  }{#1}
     }
\newbox\labelbox
\def\ugor{\goodbreak
\bigskip
\ifmargin
\hskip-1in {\box\labelbox}
\vskip-\baselineskip
\fi
\noindent\bf}
\def\skippy{\enskip\enskip}

\def\definition { \ugor Definition \sorszaminc \rm \ }   
\def\theorem {    \ugor Theorem \sorszaminc \sl \ }   
\def\prop {       \ugor Proposition \sorszaminc \sl \ }
\def\lemma {      \ugor Lemma \sorszaminc \sl \ }

\def\proof { \medskip {\noindent \it Proof.\skippy}\rm}
\def\example {    \ugor Example \sorszaminc \rm \ }

\def\nullbox{\setbox0=\null \ht0=5pt \wd0=5pt \dp0=0pt \box0}
\def\eop {\hfill \boxit{\nullbox}\goodbreak}    
\def\Eop {\vskip -\baselineskip \vskip -\belowdisplayshortskip \eop}


\def\label#1{  
\ifmargin
\vskip0in\hskip-1in {\kis#1}
\vskip-\baselineskip
\fi
\immediate
\write0{\string\redef{\string#1}{{\string\rm\the\section.\the\sorszam}}}}

\def\mark#1{  
\advance\sorszam by 1
\setbox\labelbox=\hbox{\kis #1}
\immediate
\write0{\string\redef{\string#1}{{\string\rm
\ifnum \section=0 \else\the\section.\fi
\the\sorszam}}}
\advance\sorszam by -1
}

\def\labelref#1{ 
\immediate
\write0{\string\redef{\string#1}{{\string\rm\the\sorszam}}}}

\def\cite#1{\expandafter\ifx\csname#1\endcsname\relax
    ???
    \message{!!!!!!! Missing reference '\string#1' !!!!!!!}
    \else
\csname#1\endcsname
\fi}


\def\title#1{\centerline{\nagy #1}}




\duplex


\def\mapeitoionto{\lhook\joinrel\relbar\joinrel\relbar\joinrel\longrightarrow\kern -.8em\to}
\def\mapwitoionto{\leftarrow\kern -.8em\longleftarrow\joinrel\relbar\joinrel\relbar\joinrel\rhook} 
\def\mapsonto{\Big\downarrow\lower .4ex\llap
  {$\vcenter{\hbox{$\Big\downarrow$}}$}}
\def\mapnonto{\Big\uparrow\raise .4ex\llap
  {$\vcenter{\hbox{$\Big\uparrow$}}$}}



\def\crss#1{c(#1)}

\def \slim{\mathop{\hbox{\rm s-lim}}}

%
%

\bigskip
\title{$C^{*}$-CROSSED PRODUCTS BY TWISTED}
\title{INVERSE 
SEMIGROUP 
ACTIONS}
\shorttitle{twisted semigroup actions}    

\author{N\'andor Sieben}

\bigskip

\abstract{The notions of Busby-Smith and Green type 
twisted actions are extended to discrete 
unital inverse semigroups.  The connection between the 
two types, and the connection with twisted partial 
actions, are investigated.  Decomposition theorems for 
the twisted crossed products are given.   } 


\newsection{Introduction}

\uj The work of Renault [Re1] connecting a locally 
compact groupoid to its ample inverse semigroup makes 
the study of inverse semigroups interesting.  Some of 
the results on inverse semigroups are Paterson's [Pa1] 
and Duncan and Paterson's work [DP1], [DP2] on the 
$C^{*}$-algebras of inverse semigroups, Kumjian's 
localization $C^{*}$-algebras [Kum] and Nica's $\tilde {F}$-inverse 
semigroups [Nic].  We have seen in [Sie] that the theory 
of crossed products can be generalized to inverse 
semigroups.  The strong connection between the 
$C^{*}$-algebras of locally compact groupoids and inverse 
semigroups found by Paterson [Pa2] promises a similar 
connection between the groupoid crossed products of 
[Re2] and inverse semigroup crossed products.  

Green [Gre] and Packer and Raeburn [PR] showed how to 
use twisted crossed products to decompose crossed 
product $C^{*}$-algebras.  In this paper we partially 
generalize their results to discrete inverse semigroups.  
We prove decomposition theorems for both Green and 
Busby-Smith style twisted actions.  We show that unlike 
in the group case, Green twisted actions seem slightly 
more general than Busby-Smith twisted actions.  

We show that the close connection between partial 
actions [Ex1], [McC] and inverse semigroup actions seen 
in [Sie] and [Ex3] still holds for the twisted partial 
actions of [Ex2] and Busby-Smith twisted inverse 
semigroup actions. It is a natural question to ask, whether 
there is a similar connection between Green twisted 
inverse semigroup actions and some sort of 
Green twisted partial actions.   

The research for this paper was carried out while the 
author was a student at Arizona State University under 
the supervision of John Quigg.  Part of the research 
was done during a short visit at the University of 
Newcastle.  The author is grateful to Professor Iain 
Raeburn for his hospitality.  

\newsection{Twisted inverse semigroup actions}  

\uj Twisted actions of locally compact groups were 
introduced in [BS].  The inverse semigroup 
version closely follows Exel's definition of twisted 
partial actions in [Ex2].  

Recall that a semigroup $S$ is an {\it inverse semigroup\/} if for 
every $s\in S$ there exists a unique element $s^{*}$ of $S$ so 
that $ss^{*}s=s$ and $s^{*}ss^{*}=s^{*}$.  The map $s\mapsto s^{*}$ is an 
involution.  An element $f\in S$ satisfying $f^2=f$ is called 
an {\it idempotent\/} of $S$.  The set of idempotents of an 
inverse semigroup is a semilattice. There is a natural 
partial order on $S$ defined by $s\leq t$ if and only if 
$s=ts^{*}s$.   

\definition Let $A$ be a $C^{*}$-algebra.  A {\it partial automorphism}
of $A$ is an isomorphism between two closed ideals of $A$.

\mark{taction}
\definition Let $A$ be a $C^{*}$-algebra, and let $S$ be a unital 
inverse semigroup with idempotent semilattice $E$, and 
unit $e$.  A {\it Busby-Smith twisted action }
of $S$ on $A$ is a pair $(\beta ,w)$, where for all $s\in S$, 
$\beta_s:E_{s^{*}}\to E_s$ is a partial automorphism of $A$, and for all 
$s,t\in S$, $w_{s,t}{}_{}$ is a unitary multiplier of $E_{st}$, such that 
for all $r,s,t\in S$ we have 

\item{(a)}{$E_e=A$;}
\item{(b)}{$\beta_s\beta_t=\hbox{\rm Ad }w_{s,t}\circ\beta_{st}$;} 
\item{(c)}{$ $$w_{s,t}=1_{M(E_{st})}$ if $s$ or $t$ is an idempotent; 
\item{(d)}{$\beta_r(aw_{s,t})w_{r,st}=\beta_r(a)w_{r,s}w_{rs,t}$ for all $
a\in E_{r^{*}}E_{st}$.} 

\noindent We also refer to the quadruple $(A,S,\beta ,w)$ as
a {\it Busby-Smith twisted action}.  

\uj Note that $\beta_r(aw_{s,t})$ makes sense since $a=a_1a_2$ for 
some $a_1,a_2\in E_{r^{*}}E_{st}$ and so 
$aw_{s,t}=a_1a_2w_{s,t}\in a_1E_{st}\subset E_{r^{*}}$.

Our definition is a generalization of inverse semigroup 
actions defined in [Sie].  Every inverse semigroup action 
is a trivially twisted Busby-Smith inverse semigroup 
action by taking $w_{r,s}=1_{M(E_{rs})}$.  Conversely, every 
trivially twisted Busby-Smith inverse semigroup action 
is an inverse semigroup action.  The definition is also a 
generalization of discrete twisted group actions in case 
our inverse semigroup $S$ is actually a group.  

The basic properties of Busby-Smith twisted inverse semigroup 
actions are collected in the following.  

\lemma If $(A,S,\beta ,w)$ is a Busby-Smith twisted action and 
$r,s,t\in S$, then \rm\ 
\item{(a)}{$E_s=E_{ss^{*}}$;} 
\item{(b)}{$\beta_{ss^{*}}=\hbox{\rm id}_{E_s}$;} 
\item{(c)}{$\beta_e=\hbox{\rm id}_A$;} 
\item{(d)}{$\beta_{s^{*}}=\hbox{\rm Ad }w_{s^{*},s}\circ\beta_s^{
-1}$;} 
\item{(e)}{$\beta_r(E_{r^{*}}E_s)=E_{rs}$;} 
\item{(f)}{$\beta_r(aw_{s,t}^{*})=\beta_r(a)w_{r,st}w_{rs,t}^{*}w_{
r,s}^{*}$ {\sl for all }
$a\in E_{r^{*}}E_{st}$;} 
\item{(g)}{$\beta_r(w_{s,t}a)=w_{r,s}w_{rs,t}w_{r,st}^{*}\beta_r(
a)$ {\sl for all }
$a\in E_{r^{*}}E_{st}$;} 
\item{(h)}{$\beta_r(w_{s,t}^{*}a)=w_{r,st}w_{rs,t}^{*}w_{r,s}^{*}
\beta_r(a)$ {\sl for all }
$a\in E_{r^{*}}E_{st}$;} 
\item{(i)}{$w_{s^{*}r^{*}r,s}=w_{s^{*},r^{*}rs}$;} 
\item{(j)}{$w_{s^{*},s}1_{M(E_{s^{*}r^{*}})}=w_{s^{*},r^{*}rs}$;} 
\item{(k)}{$\beta_s(w_{s^{*},s})=w_{s,s^{*}}$.} 

\label{aproperty} 

\proof\ (a) and (b) follow from the calculations 
$$E_s=\hbox{\rm dom }(\beta_s\beta_{s^{*}})=\hbox{\rm dom }(\hbox{\rm Ad }
w_{s,s^{*}}\circ\beta_{ss^{*}})=E_{ss^{*}}$$
$$\beta_f=\beta_f\beta_f\beta_f^{-1}=\hbox{\rm Ad }1_{M(E_f)}\circ
\beta_{ff}\beta_f^{-1}=\beta_f\beta_f^{-1}=\hbox{\rm id}_{E_f}\,,$$
where $f=ss^{*}$.  (c) is a special case of (b).  Using (b) we 
have (d) since 
$$\beta_{s^{*}}=\beta_{s^{*}}\beta_s\beta_s^{-1}=\hbox{\rm Ad }w_{
s^{*},s}\circ\beta_{s^{*}s}\beta_s^{-1}=\hbox{\rm Ad }w_{s^{*},s}
\circ\beta_s^{-1}.$$
We have 
$$\beta_r(E_{r^{*}}E_s)=\hbox{\rm im}\,(\beta_r\beta_s)=\hbox{\rm im}
\,(\hbox{\rm Ad }w_{r,s}\circ\beta_{rs})=E_{rs},$$
which proves (e).  Replacing $a$ by $aw_{s,t}^{*}$ in Definition 
\cite{taction}(d) gives (f).  Applying (f) we have 
(g) because 
$$\beta_r(w_{s,t}a)=\beta_r(a^{*}w_{s,t}^{*})^{*}=(\beta_r(a^{*})
w_{r,st}w_{rs,t}^{*}w_{r,s}^{*})^{*}=w_{r,s}w_{rs,t}w_{r,st}^{*}\beta_
r(a).$$
Replacing $a$ by $w_{s,t}^{*}a$ in (g) we have (h). To show (i) let 
$b\in E_{s^{*}r^{*}}$. Then by (a), (b) and (e) there is an 
$a\in E_sE_{r^{*}}=\beta_{r^{*}r}(E_{r^{*}r}E_s)=E_{r^{*}rs}$ such that 
$b=\beta_{s^{*}}(a)$. Hence
$$\eqalign{b&=\beta_{s^{*}}(aw_{r^{*}r,s})\qquad\hbox{\rm ( \cite{taction}(c) )}\cr
&=\beta_{s^{*}}(a)w_{s^{*},r^{*}r}w_{s^{*}r^{*}r,s}w_{s^{*},r^{*}
rs}^{*}\qquad\hbox{\rm ( \cite{taction}(d) )}\cr
&=bw_{s^{*}r^{*}r,s}w_{s^{*},r^{*}rs}^{*}\qquad\hbox{\rm ( (a) and  \cite{taction}(c) )}\,
,\cr}
$$
which means $w_{s^{*}r^{*}r,s}w_{s^{*},r^{*}rs}^{*}$ is the identity of $
E_{s^{*}r^{*}}$.  
This implies our statement since by (a) both $w_{s^{*}r^{*}r,s}$ 
and $w_{s^{*},r^{*}rs}$ are unitary multipliers of $E_{s^{*}r^{*}}$.  To show 
(j) let $a$ and $b$ as in the proof of part (i).  Then we 
have 
$$\eqalign{b&=\beta_{s^{*}}(aw_{s,s^{*}r^{*}rs})\qquad\hbox{\rm ( \cite{taction}(c) )}\cr
&=\beta_{s^{*}}(a)w_{s^{*},s}w_{s^{*}s,s^{*}r^{*}rs}w_{s^{*},ss^{
*}r^{*}rs}^{*}\qquad\hbox{\rm ( \cite{taction}(d) )}\cr
&=bw_{s^{*},s}1_{M(E_{s^{*}r^{*}})}w_{s^{*},r^{*}rs}^{*}\cr}
$$
and the statement follows as in part (i).  Finally, (k) 
follows from Definition \cite{taction}(c, d) if we 
extend $\beta_s$ to the multipliers of $E_{s^{*}}$.  
\eop

\uj Recall from [Pet], [CP], [How] that a {\it congruence relation }
on an inverse semigroup $S$ is an equivalence 
relation $\sim$ on $S$ such that if $s\sim t$ then $rs\sim rt$ and 
$sr\sim tr$ for all $r\in S$.  If $E$ is the idempotent semilattice 
of $S$, then the {\it kernel normal system\/} of the congruence $
\sim$ 
on $S$ is the set $\{[f]:f\in E\}$ of congruence classes 
containing idempotents.  The kernel normal system is 
exactly the set of idempotents in the quotient inverse 
semigroup $S/\!\!\sim$.  The kernel normal system determines 
the congruence relation, since 
$\sim\;=\{(s,t):ss^{*},tt^{*},st^{*}\in [f]\hbox{\rm \ for some $
f\in E$}\}$.  

If $\sim$ is an {\it idempotent-separating\/} congruence, that is, no 
two idempotents are congruent, then every equivalence 
class $[f]$ in the kernel normal system is a group with 
identity $f$.  The union $N=\cup \{[f]:f\in E\}$ is an inverse 
subsemigroup of $S$ contained in the centralizer of $E$.  It 
is also a {\it normal} subsemigroup, that is, $E\subset N$ and 
$sNs^{*}$$\subset N$ for all $s\in S$.  On the other hand, if 
$N$ is normal subsemigroup of $S$ in the centralizer of $E$ 
then $N$ determines a kernel normal system $\{[f]:f\in E\}$ of 
an idempotent-separating congruence $\sim$, where 
$[f]=\{s\in N:ss^{*}=f\}$. 
We also write $S/N$ for $S/\!\!\sim$. 
Note that a normal subsemigroup $N$ of $S$ is contained in 
the centralizer of $E$ if and only if $N$ is a {\it Clifford }
{\it semigroup}, that is, $N$ is an inverse semigroup such that 
$n^{*}n=nn^{*}$ for all $n\in N$. We call such a subsemigroup $N$ a 
{\it normal Clifford subsemigroup\/} of $S$, that is, $N$ is a 
normal subsemigroup which is also a Clifford 
semigroup. We thus get a bijective correspondence 
between idempotent-separating congruences on $S$ and 
normal Clifford subsemigroups of $S$.
In the theory of twisted group actions, normal subgroups 
play an important role.  For inverse semigroups the 
situation is more complicated, but normal Clifford 
subsemigroups give an appropriate substitute for normal 
subgroups.  

\definition If $N$ is a normal Clifford subsemigroup of the 
inverse semigroup $S$ with idempotent semilattice $E$, then 
a cross-section $c:S/N\to S$ is called {\it order-preserving\/} if 
$c([f])=f$ for all $f\in E$ and $[s]\leq [t]$ implies $c([s])\leq 
c([t])$ 
for all $s,t\in S$.  

\prop Let $T$ be an inverse semigroup with idempotent 
semilattice $E$, and let $N$ be a normal Clifford 
subsemigroup of $T$.  Let $S=T/N$, and suppose there is 
an order-preserving cross-section $c:S\to T$.  For $s\in S$ 
define 
$$E_s=\overline {\hbox{\rm span}}\,\cup \{[f]:f\leq\crss s\crss s^{
*},f\in E\},$$
where the closed linear span is taken in $C^{*}(N)$.  
Then each $E_s$ is a closed ideal of $C^{*}(N)$.  For $r,s\in S$ 
define 
$$\eqalign{&\beta_s:E_{s^{*}}\to E_s\ \ {\rm b}{\rm y}\ \ \beta_s
(a)=\crss sa\crss s^{*},\quad\hbox{\rm and}\cr
&w_{r,s}=\crss r\crss s\crss{rs}^{*}\in M(E_{rs}).\cr}
$$
Then $(C^{*}(N),S,\beta ,w)$ is a Busby-Smith twisted action.  

\label{tactionexample}

\proof First notice that since $\crss{rs}^{*}\crss{rs}$ is the 
identity of the congruence class containing 
$(c(r)c(s))^{*}c(r)c(s)$, 
we have 
$\crss r\crss s=\crss r\crss s\crss{rs}^{*}\crss{rs}=w_{r,s}\crss{
rs}$ 
for all $r,s\in S$.  Also $\crss s\crss s^{*}=\crss{ss^{*}}$ because 
$c([f])=f$ for all $f\in E$, 
and similarly 
$\crss s^{*}\crss s=\crss{s^{*}s}$.  Each $E_s$ is a right ideal because 
if $a\in N\cap E_s$, that is, $a\in [f]$ for some $f\leq\crss s\crss 
s^{*}$ and 
$b\in [g]$, then $ab\in [fg]\subset E_s$ since $fg\leq\crss s\crss 
s^{*}$. 
A similar argument shows that $E_s$ is also a left ideal.  It 
is clear that each $\beta_s$ is an isomorphism and each $w_{r,s}$ is 
a unitary multiplier of $E_{rs}$.  It remains to 
check the conditions of Definition \cite{taction}.  
Condition (a) holds since $E_e=C^{*}(N)$.  

To check (b), let $a\in\hbox{\rm dom}\,(\beta_r\beta_s)\cap\hbox{\rm dom}
\,(\beta_{rs})$. Then 
$$\beta_r\beta_s(a)=\crss r\crss sa\crss s^{*}\crss r^{*}=w_{r,s}\crss{
rs}a\crss{rs}^{*}w_{r,s}^{*}=\hbox{\rm Ad }w_{r,s}\circ\beta_{rs}
(a).$$
We need to show that $\beta_r\beta_s$ and $\hbox{\rm Ad }w_{r,s}\circ
\beta_{rs}$ have the same domain.  
First we show that ~$\hbox{\rm dom}\,\beta_r\beta_s\subset\hbox{\rm dom}\,
\beta_{rs}$.  If $a\in\hbox{\rm dom}\,\beta_r\beta_s$, then 
$a=\lim_ia_i$ for some $a_i\in\hbox{\rm span}\,\cup \{[f]:f\leq c
(s)^{*}c(s)\}$ and $c(s)ac(s)^{*}\in E_{r^{*}}$.  
Hence $c(s)ac(s)^{*}=\lim_jb_j$ for some $b_j\in\hbox{\rm span}\,
\cup \{[f]:f\leq c(r)^{*}c(r)\}$.  Since 
$$c(s)^{*}c(s)ac(s)^{*}c(s)=\lim_ic(s)^{*}c(s)a_ic(s)^{*}c(s)=\lim_
ia_i=a\,,$$
we have $a=\lim_jc(s)^{*}b_jc(s)$.  It suffices to show that 
$c(s)^{*}b_jc(s)\in\hbox{\rm span}\,\cup \{[f]:f\leq c(rs)^{*}c(r
s)\}$.  This follows from the fact 
that if $b\in [f]$ such that $f\leq c(r)^{*}c(r)$, then $c(s)^{*}
bc(s)\in [c(s)^{*}fc(s)]$ and 
$c(s)^{*}fc(s)\leq c(s)^{*}c(r)^{*}c(r)c(s)=c(rs)^{*}c(rs)$.  

Next we show $\hbox{\rm dom}\,\beta_r\beta_s\supset\hbox{\rm dom}\,
\beta_{rs}$.  If $a\in\hbox{\rm dom}\,\beta_{rs}$, then $a=\lim_i
a_i$ for 
some $a_i\in\hbox{\rm span}\,\cup \{[f]:f\leq c(rs)^{*}c(rs)\}$.  Since 
$$c(rs)^{*}c(rs)=c(s)^{*}c(r)^{*}c(r)c(s)\leq c(s)^{*}c(s)\,,$$
$a_i\in\hbox{\rm span}\,\cup \{[f]:f\leq c(s)^{*}c(s)\}$ and so $
a\in E_{s^{*}}$.  We have 
$c(s)ac(s)^{*}=c(s)\lim_ia_ic(s)^{*}=\lim_ic(s)a_ic(s)^{*}$, so it remains to check 
that $c(s)a_ic(s)^{*}\in\cup \{[f]:f\leq c(r)^{*}c(r)\}.$ This is true since if $
b\in [f]$ such 
that $f\leq c(rs)^{*}c(rs)$, then $c(s)bc(s)^{*}\in [c(s)fc(s)^{*}
]$ and 
$$c(s)fc(s)^{*}\leq c(s)c(rs)^{*}c(rs)c(s)^{*}=c(s)c(s)^{*}c(r)^{
*}c(r)c(s)c(s)^{*}\leq c(r)^{*}c(r)\,.$$

To check (c), fix $a\in N\cap E_{[f]s}$ for some $f\in E$.  Then 
$a\in [g]$ for some idempotent $g\in T$ such that 
$$g\leq\crss{[f]s}\crss{[f]s}^{*}=c([f]ss^{*}[f])=c([f]ss^{*})=f\crss{
ss^{*}}.$$
Since $c$ is 
order-preserving and $[f]s\leq s$, we have 
$$w_{[f],s}=\crss{[f]}\crss s\crss{[f]s}^{*}\ge f\crss{[f]s}\crss{
[f]s}^{*}=f\crss{ss^{*}}$$
and so $aw_{[f],s}=w_{[f],s}a=a$.  The other part of (c) 
follows similarly.  

To check (d), fix $a\in N\cap { E}_{r^{*}}{ E}_{st}$.  Since the classes in the 
kernel normal system are disjoint, $a\in [f]$ for some $f\in T$ 
such that 
$f\leq\crss{r^{*}}\crss{r^{*}}^{*}=\crss{r^{*}r}=\crss r^{*}\crss 
r$ and 
$f\leq\crss{st}\crss{st}^{*}$.  Now we have 
$$\eqalign{\crss raw_{s,t}\crss r^{*}w_{r,st}&=\crss ra\crss s\crss 
t\crss{st}^{*}\crss r^{*}\crss r\crss{st}\crss{rst}^{*}\,.\cr}
$$
Since $a\crss s\crss t\crss{st}^{*}\in [f]st(st)^{*}$, the domain 
projection of $a\crss s\crss t\crss{st}^{*}$ is 
$f\crss{st}\crss{st}^{*}=f$.  Hence 
$a\crss s\crss t\crss{st}^{*}\crss r^{*}\crss r=a\crss s\crss t\crss{
st}^{*}$, 
so 
$$\eqalign{\crss raw_{s,t}\crss r^{*}w_{r,st}&=\crss ra\crss s\crss 
t\crss{st}^{*}\crss{st}\crss{rst}^{*}=\crss ra\crss s\crss t\crss{
rst}^{*}\cr
&=\crss ra\crss r^{*}\crss r\crss s\crss t\crss{rst}^{*}\cr
&=\crss ra\crss r^{*}\crss r\crss s\crss{rs}^{*}\crss{rs}\crss t\crss{
rst}^{*}\cr
&=\crss ra\crss r^{*}w_{r,s}w_{rs,t}.\cr}
$$
\Eop

\uj In Proposition \cite{tactionexproof} we further 
investigate the action in the previous lemma.  

\example If in Proposition \cite{tactionexample} $N$ is the 
idempotent semilattice of $T$, then the (trivially) twisted 
Busby-Smith
action becomes the {\it canonical action\/} of the inverse 
semigroup on its semilattice [Sie].  \label{canonicalaction} 

\uj Proposition \cite{tactionexample} shows the 
significance of an order-preserving cross-section.  
Unfortunately, these cross-sections do not always exist, 
as we can see in the following.  

\example Let $s=$$(a_i,\cdots, a_n)$ denote the partial bijection in 
the symmetric inverse semigroup ${\cal I}(\{1,\dots,n\})$ whose 
domain is $\{i:a_i\ne 0\}$ and $s(i)=a_i$ for all $i\in\hbox{\rm dom }
s$.  Let 
$n=6$ and consider the 19-element inverse subsemigroup 
generated by the elements 
$$r=(1,4,5,0,0,0),\quad s=(0,5,4,0,0,6)\,.$$
It is tedious but not hard to check that we can define a kernel 
normal system consisting of $\{s^{*}r,s^{*}rs^{*}r\}$, $\{rs^{*},
rs^{*}rs^{*}\}$ and the 
singleton sets containing idempotents other than $s^{*}rs^{*}r$ and 
$rs^{*}rs^{*}$, and the congruence determined by this kernel 
normal system separates idempotents. Then 
$[ss^{*}r]\leq [r]=\{r\}$ 
and $[ss^{*}r]\leq [s]=\{s\}$, but we cannot choose a representative of 
$[ss^{*}r]$ which is less than both $s$ and $r$.

\uj The situation is not as bad as it seems.  In many 
cases we have an order-preserving cross section.  Recall 
[Nic] that an $\widetilde {F}${\it -inverse semigroup\/} is a unital inverse 
semigroup in which every non-zero element is majorized 
by a unique maximal element.  

\prop Let $S$ be a unital inverse semigroup with identity 
$e$, and let $\sim$ be an idempotent-separating congruence on 
$S$ such that $T=S/\!\!\sim$ is an $\widetilde{F}$-inverse 
semigroup.  Denote the set of maximal elements in $T$ by 
${\cal M}$.  Choose a cross-section $c:{\cal M}\to S$ with $c([e]
)=e$.  
Extend $c$ to $T$ by defining $c(0)=0$ and $c(t)=c(m_t)f$ for 
$t\neq 0$, where $m_t\in {\cal M}$ is the unique maximal element 
majorizing $t$ and $f$ is the idempotent in $[t]^{*}[t]$.  Then $
c$ 
is an order-preserving cross-section.  

\proof It is clear that for all non-zero idempotents $f\in S$ 
we have $m_{[f]}=[e]$ and hence $c([f])=f$.  If $0\neq s\leq t$ in $
S/\!\!\sim$ 
then $m_s=m_t$ and so 
$$\eqalign{c(t)c(s)^{*}c(s)&=c(m_t)c(t^{*}t)c(s^{*}s)c(m_s)^{*}c(
m_s)c(s^{*}s)\cr
&=c(m_s)c(t^{*}t)c(s^{*}s)=c(m_s)c(s^{*}s)=c(s)\,.\cr}
$$
\Eop

\newsection{Busby-Smith Twisted crossed products}

\uj
Given a Busby-Smith twisted inverse semigroup action 
$(A,S,\beta ,w)$, define multiplication and involution on 
the closed subspace
$$L=\{x\in l^1(S,A):x(s)\in E_s\hbox{\rm \ for all }s\in S\}$$
of $l^1(S,A)$
by 
$$\eqalign{\big(x*y\big)(s)&=\sum_{rt=s}\beta_r(\beta_r^{-1}(x(r)
)y(t))w_{r,t}\quad\hbox{\rm and}\quad x^{*}(s)=w_{s,s^{*}}\beta_s
(x(s^{*})^{*}).\cr}
$$
We will show that the multiplication is well-defined in 
the proof of Proposition \cite{banachref} below. 
For $a\in E_s$ we are going to denote by $a\delta_s$ the function in $
L$ taking 
the value $a$ at $s$ and zero at every other element of $S$, 
so that $L$ is the closed span of the $a\delta_s$.  
Then we have 
$$\eqalign{a_s\delta_s*a_t\delta_t&=\beta_s(\beta_s^{-1}(a_s)a_t)
w_{s,t}\delta_{st}\cr
(a\delta_s)^{*}&=\beta_s^{-1}(a^{*})w_{s^{*},s}^{*}\delta_{s^{*}}
.\cr}
$$
Notice that $a_e\delta_e*a_t\delta_t=a_ea_t\delta_t$ and 
$a_s\delta_s*a_e\delta_e=\beta_s(\beta_s^{-1}(a_s)a_e)\delta_s$.
 
\prop$L$ is a Banach $^{*}$-algebra.  

\label{banachref}

\proof First notice that $\big(x*y\big)(s)\in E_s$ by Proposition 
\cite{aproperty}(e) and so to show that the product is 
well-defined we only need to check that $\|x*y\|$ is finite:
$$\eqalign{\|x*y\|&=\sum_{s\in S}\sum_{rt=s}\|\beta_r(\beta_r^{-1}
(x(r))y(t))w_{r,t}\|\cr
&\leq\sum_{s\in S}\sum_{rt=s}\|\beta_r^{-1}(x(r))y(t)\|\leq\sum_{
r\in S}\sum_{t\in S}\|\beta_r^{-1}(x(r))\|\|y(t)\|\cr
&\leq\sum_{r\in S}\|\beta_r^{-1}(x(r))\|\sum_{t\in S}\|y(t)\|\leq
\|x\|\|y\|.\cr}
$$
One can easily check that $\|x^{*}\|=\|x\|$ and so $x^{*}\in L$. For 
$a\in E_s$ we have  
$$\eqalign{(a\delta_s)^{**}&=\beta_{s^{*}}^{-1}(w_{s^{*},s}\beta_
s^{-1}(a))w_{s,s^{*}}^{*}\delta_s\cr
&=w_{s,s^{*}}^{*}\beta_s((w_{s^{*},s}^{*}\beta_{s^{*}}(x(s)^{*}))^{
*})\cr
&=w_{s,s^{*}}^{*}\beta_s(w_{s^{*},s}\beta_s^{-1}(a))\delta_s\qquad\hbox{\rm (Proposition \cite{aproperty}(d) )}\cr
&=w_{s,s^{*}}^{*}w_{s,s^{*}}w_{ss^{*},s}w_{s,s^{*}s}^{*}\beta_s(\beta_
s^{-1}(a))\delta_s\qquad\hbox{\rm (Proposition \cite{aproperty}(g) )}\cr
&=a\delta_s.\cr}
$$
Next we show that the multiplication is associative.  Let 
$a$, $b$ and $c$ be elements of $E_r$, $E_s$ and $E_t$, respectively, 
and let $u_{\lambda}$ be an approximate identity for $E_{s^{*}}$.  Then we 
have 
$$\eqalign{(a\delta_r*b\delta_s)*c\delta_t&=\beta_{rs}(\beta_{rs}^{
-1}(\beta_r(\beta_r^{-1}(a)b)w_{r,s})c)w_{rs,t}\delta_{rst}\cr
&=\lim_{\lambda}\beta_r(\beta_r^{-1}(a)b)w_{r,s}\beta_{rs}(u_{\lambda}
c)w_{rs,t}\delta_{rst}\cr
&=\lim_{\lambda}\beta_r(\beta_r^{-1}(a)b)\beta_r(\beta_s(u_{\lambda}
c))w_{r,s}w_{rs,t}\delta_{rst}\cr
&=\lim_{\lambda}\beta_r(\beta_r^{-1}(a)b\beta_s(u_{\lambda}c))w_{
r,s}w_{rs,t}\delta_{rst}\cr
&=\lim_{\lambda}\beta_r(\beta_r^{-1}(a)\beta_s(\beta_s^{-1}(b)u_{
\lambda}c))w_{r,s}w_{rs,t}\delta_{rst}\,\cr
&=\beta_r(\beta_r^{-1}(a)\beta_s(\beta_s^{-1}(b)c))w_{r,s}w_{rs,t}
\delta_{rst}\cr
&=\beta_r(\beta_r^{-1}(a)\beta_s(\beta_s^{-1}(b)c)w_{s,t})w_{r,st}
\delta_{rst}\qquad\hbox{\rm ( \cite{taction}(c) ),}\cr
&=a\delta_r*(b\delta_s*c\delta_t)\,.\cr}
$$
Finally we show that the involution is 
anti-multiplicative.  Fix $a\in E_r$ and $b\in E_s$.  Then we 
have 
$$\eqalign{&(a\delta_r*b\delta_s)^{*}\cr
&=(\beta_r(\beta_r^{-1}(a)b)w_{r,s}\delta_{rs})^{*}\cr
&=\beta_{rs}^{-1}(w_{r,s}^{*}\beta_r(b^{*}\beta_r^{-1}(a^{*})))w_{
s^{*}r^{*},rs}^{*}\delta_{s^{*}r^{*}}\cr
&=w_{s^{*}r^{*},rs}^{*}\beta_{s^{*}r^{*}}(w_{r,s}^{*}\beta_r(b^{*}
\beta_r^{-1}(a^{*})))\delta_{s^{*}r^{*}}\qquad\hbox{\rm ( \cite{aproperty}(d) )}\cr
&=w_{s^{*}r^{*},rs}^{*}\beta_{s^{*}r^{*}}(w_{r,s}^{*}\beta_r(\beta_
s(\beta_s^{-1}(b^{*}\beta_r^{-1}(a^{*})))))\delta_{s^{*}r^{*}}\cr
&=w_{s^{*}r^{*},rs}^{*}\beta_{s^{*}r^{*}}(w_{r,s}^{*}w_{r,s}\beta_{
rs}(\beta_s^{-1}(b^{*}\beta_r^{-1}(a^{*})))w_{r,s}^{*})\delta_{s^{
*}r^{*}}\qquad\hbox{\rm ( \cite{taction}(b) )}\cr
&=w_{s^{*}r^{*},rs}^{*}\beta_{s^{*}r^{*}}(\beta_{rs}(\beta_s^{-1}
(b^{*}\beta_r^{-1}(a^{*}))))w_{s^{*}r^{*},rs}w_{s^{*}r^{*}r,s}^{*}
w_{s^{*}r^{*},r}^{*}\delta_{s^{*}r^{*}}\qquad\hbox{\rm ( \cite{aproperty}(f) )}\cr
&=\beta_{rs}^{-1}(\beta_{rs}(\beta_s^{-1}(b^{*}\beta_r^{-1}(a^{*}
))))w_{s^{*}r^{*}r,s}^{*}w_{s^{*}r^{*},r}^{*}\delta_{s^{*}r^{*}}\,
.\cr}
$$
On the other hand, we have
$$\eqalign{&(b\delta_s)^{*}*(a\delta_r)^{*}\cr
&=(\beta_s^{-1}(b^{*})w_{s^{*},s}^{*}\delta_{s^{*}})*(\beta_r^{-1}
(a^{*})w_{r^{*},r}^{*}\delta_{r^{*}})\cr
&=\beta_{s^{*}}(\beta_{s^{*}}^{-1}(\beta_s^{-1}(b^{*})w_{s^{*},s}^{
*})\beta_r^{-1}(a^{*})w_{r^{*},r}^{*})w_{s^{*},r^{*}}\delta_{s^{*}
r^{*}}\cr
&=\beta_{s^{*}}(w_{s,s^{*}}^{*}\beta_s(\beta_s^{-1}(b^{*})w_{s^{*}
,s}^{*})w_{s,s^{*}}\beta_r^{-1}(a^{*})w_{r^{*},r}^{*})w_{s^{*},r^{
*}}\delta_{s^{*}r^{*}}\cr
&=\beta_{s^{*}}(w_{s,s^{*}}^{*}\beta_s(\beta_s^{-1}(b^{*}))w_{s,s^{
*}s}w_{ss^{*},s}^{*}w_{s,s^{*}}^{*}w_{s,s^{*}}\beta_r^{-1}(a^{*})
w_{r^{*},r}^{*})w_{s^{*},r^{*}}\delta_{s^{*}r^{*}}\quad\hbox{\rm ( \cite{aproperty}(f) )}\cr
&=\beta_{s^{*}}(w_{s,s^{*}}^{*}b^{*}\beta_r^{-1}(a^{*})w_{r^{*},r}^{
*})w_{s^{*},r^{*}}\delta_{s^{*}r^{*}}\qquad\hbox{\rm ( \cite{taction}(c) )}\cr
&=w_{s^{*},ss^{*}}w_{s^{*}s,s^{*}}^{*}w_{s^{*},s}^{*}\beta_{s^{*}}
(b^{*}\beta_r^{-1}(a^{*})w_{r^{*},r}^{*})w_{s^{*},r^{*}}\delta_{s^{
*}r^{*}}\qquad\hbox{\rm ( \cite{aproperty}(h) )}\cr
&=w_{s^{*},s}^{*}\beta_{s^{*}}(b^{*}\beta_r^{-1}(a^{*}))w_{s^{*},
r^{*}r}w_{s^{*}r^{*},r}^{*}w_{s^{*},r^{*}}^{*}w_{s^{*},r^{*}}\delta_{
s^{*}r^{*}}\qquad\hbox{\rm ( \cite{aproperty}(f) )}\cr
&=\beta_s^{-1}(b^{*}\beta_r^{-1}(a^{*}))w_{s^{*},s}^{*}1_{M(E_{s^{
*}r^{*}r})}w_{s^{*}r^{*},r}^{*}\delta_{s^{*}r^{*}}\qquad\hbox{\rm ( \cite{aproperty}(d) ),}\cr}
$$
and so the equality $(a\delta_r*b\delta_s)^{*}=(b\delta_s)^{*}*(a
\delta_r)^{*}$ follows 
from Proposition \cite{aproperty}(i) and (j).  \eop

\definition A {\it covariant representation\/} of a Busby-Smith twisted 
action $(A,S,\beta ,w)$ is a triple $(\pi ,v,H)$, where $\pi$ 
is a nondegenerate representation of $A$ on the Hilbert 
space $H$ and $v_s$ is a partial isometry for all $s\in S$, such 
that for all $r,s\in S$ we have 
\item{(a)}{$\pi (\beta_s(a))=v_s\pi (a)v_s^{*}\quad$for $a\in E_{
s^{*}}$;} 
\item{(b)}{$v_rv_s=\pi (w_{r,s})v_{rs}$;} 
\item{(c)}{$v_s$ has initial space $\pi (E_{s^{*}})H$ and final space 
$\pi (E_s)H$.} 

\noindent To evaluate $\pi (w_{r,s})$ we extend $\pi$ to the 
enveloping von Neumann algebra $A^{**}$ of $A$.  

\label{covrepdef}

\uj We sometimes use the shortened notation $(\pi,v)$ if we do not need
a symbol for the Hilbert space. 
Note that the above definition is a generalization of a 
covariant representation of an inverse semigroup 
action [Sie], which is a trivially twisted 
Busby-Smith twisted action.  

\def \detail#1{} 

\detail{ Also note that $\pi (E_s)H$ is closed in $H$.  This is a 
consequence of the Cohen-Hewitt factorization theorem, 
which we state in the following form.  
\theorem(Cohen-Hewitt) If $A$ is a Banach algebra with 
bounded approximate identity, $V$ is a Banach space and 
$\pi :A\to B(V)$ is a bounded homomorphism then 
$\pi (A)V=\{\pi (a)v:a\in A,v\in V\}$ is a closed subspace of $V$. 
\lemma If $I$ is a closed ideal of the $C^{*}$-algebra $A$, $m$ is a 
unitary multiplier of $I$ and $\pi :A\to B(H)$ is a representation 
of $A$ then $\pi (m)$ is a partial isometry with initial and 
final space $\pi (I)H$.  
\label{multlemma}
\proof$\pi (m)$ is a partial isometry since 
$\pi (m)\pi (m)^{*}\pi (m)=\pi (m)$.  If $a\in I$ and $h\in H$ then 
$\pi (a)h=\pi (mm^{-1}a)h=\pi (m)\pi (m^{-1}a)h\in\pi (m)H$.  On the other 
hand, $\pi (m)$ is in the strong operator closure of $\pi (I)$ and 
so there is a net $\{a_{\lambda}\}\subset I$ so that $\pi (a_{\lambda}
)h\to\pi (m)h$ for all 
$h\in H$.  Thus $\pi (m)H\subset\pi (I)H$ and so $\pi (m)$ has the required 
final space. The statement about the initial space 
follows from the fact that $\pi (m)^{*}=\pi (m^{-1})$.   
}

\prop Let $(\pi ,v,H)$ be a covariant representation. If
$s\in S$ and $f$ is an idempotent in $S$ then  
\rm
\item{(a)}{$v_f$ is the orthogonal projection onto $\pi (E_f)H$;} 
\item{(b)}{$v_e=1_{B(H)}$;} 
\item{(c)}{$v_{s^{*}}=\pi (w_{s^{*},s})v_s^{*}$;} 
\item{(d)}{$v_s^{*}=v_{s^{*}}\pi (w_{s,s^{*}}^{*})$;} 
\item{(e)}{$v_s^{*}=\pi (w_{s^{*},s}^{*})v_{s^{*}}$.} 

\label{vproperty}

\proof We have 
$v_f=v_fv_fv_f^{*}=\pi (w_{f,f})v_fv_f^{*}=\pi (1_{M(E_f)})$.  Since 
$E_e=A$, (b) is a special case of (a).  (c) follows from the 
calculation $v_{s^{*}}=v_{s^{*}}v_sv_s^{*}=\pi (w_{s^{*},s})v_{s^{
*}s}v_s^{*}=\pi (w_{s^{*},s})v_s^{*}$.  
We can get (d) from (c) upon taking adjoints.  Finally using 
(c) we have 
$\pi (w_{s^{*},s}^{*})v_{s^{*}}=\pi (w_{s^{*},s}^{*}w_{s^{*},s})v_
s^{*}=\pi (1_{E_{s^{*}s}})v_s^{*}=\pi (1_{E_{s^{*}}})v_s^{*}=v_s^{
*}$.  
\eop

\definition Let $(\pi ,v,H)$ be a covariant representation.  
Define $\pi\times v:L\to B(H)$ by 
$$\big(\pi\times v\big)(x)=\sum_{s\in S}\pi (x(s))v_s\,,$$
where the series converges in norm.

\prop$\pi\times v$ is a nondegenerate representation of $L$.
\label{hivcpv}

\proof It is clear that $\pi\times v$ is linear.  It suffices 
to show multiplicativity for $a\delta_s$ and $b\delta_t$ where $a
\in E_s$ 
and $b\in E_t$:
$$\eqalign{\big(\pi\times v\big)(a\delta_s*b\delta_t)&=\big(\pi\times 
v\big)(\beta_s(\beta_s^{-1}(a)b)w_{s,t}\delta_{st})\cr
&=\pi (\beta_s(\beta_s^{-1}(a)b)w_{s,t})v_{st}\cr
&=v_s\pi (\beta_s^{-1}(a)b)v_s^{*}\pi (w_{s,t})v_{st}\cr
&=v_s\pi (w_{s^{*},s}^{*}\beta_{s^{*}}(a)w_{s^{*},s})\pi (b)v_s^{
*}v_sv_t\cr
&=v_s\pi (w_{s^{*},s}^{*})v_{s^{*}}\pi (a)v_{s^{*}}^{*}\pi (w_{s^{
*},s})\pi (b)v_s^{*}v_sv_t\cr
&=v_sv_s^{*}\pi (a)v_s\pi (w_{s^{*},s}^{*})\pi (w_{s^{*},s})\pi (
b)v_s^{*}v_sv_t\cr
&=\pi (a)v_s\pi (b)v_t.\cr}
$$
To show that $\pi\times v$ preserves adjoints let $a\in E_s$. Then 
we have
$$\eqalign{\big(\pi\times v\big)(a\delta_s)^{*}&=(\pi (a)v_s)^{*}
=v_s^{*}\pi (a^{*})=v_s^{*}v_s\pi (\beta_s^{-1}(a^{*}))v_s^{*}\cr
&=\pi (\beta_s^{-1}(a^{*}))v_s^{*}=\pi (\beta_s^{-1}(a^{*}))\pi (
w_{s^{*},s}^{*})\pi (w_{s^{*},s})v_s^{*}\cr
&=\pi (\beta_s^{-1}(a^{*})w_{s^{*},s}^{*})v_{s^{*}}=\big(\pi\times 
v\big)(\beta_s^{-1}(a^{*})w_{s^{*},s}^{*}\delta_{s^{*}})\cr
&=\big(\pi\times v\big)((a\delta_s)^{*}).\cr}
$$
If $\{u_{\lambda}\}$ is a bounded 
approximate identity for $A$ then 
$\{u_{\lambda}\delta_e\}$ is a bounded approximate identity for $
L$, since for  
$a\in E_s$ we have 
$\lim_{\lambda}u_{\lambda}\delta_e*a\delta_s=\lim_{\lambda}u_{\lambda}
a\delta_s=a\delta_s$, and 
$\lim_{\lambda}a\delta_s*u_{\lambda}\delta_e=\lim_{\lambda}\beta_
s(\beta_s^{-1}(a)u_{\lambda})\delta_s=\beta_s(\beta_s^{-1}(a))\delta_
s=a\delta_s$. 
Since $\pi$ is a nondegenerate representation,
$\big(\pi\times v\big)(u_{\lambda}\delta_e)=\pi (u_{\lambda})$ converges strongly to $
1_{B(H)}$ 
and so $(\pi\times v)$ is nondegenerate.  
\eop
 
\definition Let $(A,S,\beta ,w)$ be a Busby-Smith twisted action.  
Define a $C^{*}$-seminorm $\|\cdot\|_c$ on $L$ by 
$$\|x\|_c=\sup\{\|(\pi\times v)(x)\|:(\pi ,v)\hbox{\rm \ is a covariant representation of }
(A,S,\beta ,w)\}.$$
Let $I=\{x\in L:\|x\|_c=0\}$.
The {\it Busby-Smith twisted crossed product }
$A\times_{\beta ,w}S$ is the 
completion of the quotient $L/I$ with respect to $\|.\|_c$.  
We denote the quotient map of $L$ onto $L/I$ by $\Phi$.  

\label {crospdef}

\uj Note that for trivially twisted Busby-Smith actions this 
definition gives the inverse semigroup crossed product of 
[Sie].

\uj Let $(\pi ,v)$ be a covariant representation of $(A,S,\beta ,
w),$ and 
let $\pi\times v$ be the associated representation of $L$.  Since 
$\hbox{\rm ker }\Phi\subset\hbox{\rm ker }\pi\times v$, we can factor $
\pi\times v$ through the 
quotient $L/I$ and extend to $A\times_{\beta ,w}S$ by continuity.  We 
denote this extension also by $\pi\times v$.  
Thus every covariant representation gives a 
nondegenerate representation of the crossed product.  
Proposition \cite{slimref} will show how to reverse this 
process.  The following lemma shows that the ideal $I$ 
may be nontrivial{\it .  }

\lemma If $s\leq t$ in $S$, that is, $s=ft$ for some 
idempotent $f\in S$, then $\Phi (a\delta_s)=\Phi (a\delta_t)$ for all $
a\in E_s$.  In 
particular $\Phi (a\delta_s)=\Phi (a\delta_e)$ if $s$ is an idempotent.  
\label {hiv2} 

\proof\ It is clear that $a\in E_t$.  If $(\pi ,v)$ is a covariant 
representation of $(A,S,\beta ,w)$ then 
$$\eqalign{\big(\pi\times v\big)(a\delta_{ft})&=\pi (a)v_{ft}=\pi 
(a)\pi (w_{f,t}^{*})v_fv_t\cr
&=\pi (a)v_t\qquad\hbox{\rm (Proposition \cite{vproperty}(a) )}\cr
&=\big(\pi\times v\big)(a\delta_t),\cr}
$$
which shows that $\Phi (a\delta_s-a\delta_t)=0.$ The second statement 
follows from the fact that $s=se$.  \eop

\uj In spite of the above lemma, we identify $a\delta_s$ with its 
image in $A\times_{\beta ,w}S$.

\prop Let $(\Pi ,H)$ be a nondegenerate representation of 
$A\times_{\beta ,w}S$.  Define a representation $\pi$ of $A$ on $
H$ and a map
$v:S\to B(H)$ by 
$$\eqalign{\pi (a)&=\Pi (a\delta_e)\quad\hbox{\rm and}\quad v_s=\slim\limits_{
\lambda}\Pi (u_{\lambda}\delta_s)\,,\cr}
$$
where $\{u_{\lambda}\}$ is an approximate identity for $E_s$ and $\slim$ 
denotes strong operator limit.  Then 
$(\pi ,v,H)$ is a covariant representation of $(A,S,\beta ,w)$.  
\label{slimref}

\proof$\pi$ is a nondegenerate representation, since $\{u_{\lambda}
\delta_e\}$ is 
an approximate identity for $A\times_{\beta ,w}S$ whenever $\{u_{
\lambda}\}$ is an 
approximate identity for $A$.
We show that $v_s$ is well-defined.  If 
$h\in\pi (E_{s^{*}})H$ then $h=\Pi (a\delta_e)k$ for some $a\in E_{
s^{*}}$ and $k\in H$.  
Hence
$$\eqalign{\lim_{\lambda}\Pi (u_{\lambda}\delta_s)h&=\lim_{\lambda}
\Pi (u_{\lambda}\delta_s)\Pi (a\delta_e)k=\lim_{\lambda}\Pi (u_{\lambda}
\delta_s*a\delta_e)k\cr
&=\lim_{\lambda}\Pi (\beta_s(\beta_s^{-1}(u_{\lambda})a)\delta_s)
k=\Pi (\beta_s(a)\delta_s)k,\cr}
$$
since $\beta_{s^{*}}(u_{\lambda})$ is an approximate 
identity for $E_{s^{*}}$.  Note that the limit is independent of 
the choice of $\{u_{\lambda}\}$ since the expression $h=\Pi (a\delta_
e)k$ was.  
On the other hand if $h\perp\pi (E_{s^{*}})H$ then 
$$\eqalign{\lim_{\lambda}\Pi (u_{\lambda}\delta_s)h&=\lim_{\lambda}
\Pi (\beta_s(\beta_s^{-1}(\sqrt {u_{\lambda}})\beta_s^{-1}(\sqrt {
u_{\lambda}}))\delta_s)h\cr
&=\lim_{\lambda}\Pi (\sqrt {u_{\lambda}}\delta_s*\beta_s^{-1}(\sqrt {
u_{\lambda}})\delta_e)h\cr
&=\lim_{\lambda}\Pi (\sqrt {u_{\lambda}}\delta_s)\pi (\beta_s^{-1}
(\sqrt {u_{\lambda}}))(h).\cr}
$$
But $\pi (E_{s^{*}})h=0$, so 
$v_s(h)=0$.  Hence $v_s$ is well-defined. Clearly 
$v_s$ is a bounded linear 
transformation, and if $f$ is an idempotent then $v_f$ is the 
orthogonal projection onto $\pi (E_f)H$. 
Notice that  
$$\eqalign{v_s^{*}&=\slim_{\lambda}\Pi (u_{\lambda}\delta_s)^{*}=\slim_{
\lambda}\Pi (\beta_s^{-1}(u_{\lambda})w_{s^{*},s}^{*}\delta_{s^{*}}
)=\pi (w_{s^{*},s}^{*})v_{s^{*}}.\cr}
$$
For $s$, $t\in S$ let $\{u_{\lambda}^s\}$ and $\{u_{\mu}^t\}$ be bounded approximate 
identities for $E_s$ and $E_t$ respectively. Then
$$\eqalign{v_sv_t&=\slim_{\lambda ,\mu}\Pi (u_{\lambda}^s\delta_s*
u_{\mu}^t\delta_t)=\slim\Pi (\beta_s(\beta_{s^{*}}(u_{\lambda}^s)u_{
\mu}^t)w_{s,t}\delta_{st})\cr
&=\Pi (w_{s,t}\delta_{st})=\slim\Pi (w_{s,t}\beta_s(\beta_{s^{*}}(
u_{\lambda}^s)u_{\mu}^t)\delta_{st})\cr
&=\slim\Pi (w_{s,t}\delta_e)\Pi (\beta_s(\beta_{s^{*}}(u_{\lambda}^
s)u_{\mu}^t)\delta_{st})=\pi (w_{s,t})v_{st}\,,\cr}
$$
since the net $\{\beta_s(\beta_{s^{*}}(u_{\lambda}^s)u_{\mu}^t)\}$ with the product direction 
is an approximate identity for $\beta_s(E_{s^{*}}E_t)=E_{st}$ (using 
boundedness of $\{u_{\lambda}^s\}$ and $\{u_{\mu}^t\}$). Thus $v$ is multiplicative.
We have $v_s^{*}v_s=\pi (w_{s^{*},s}^{*})v_{s^{*}}v_s=v_{s^{*}s}$, which is the projection 
onto $\pi (E_{s^{*}s})H=\pi (E_{s^{*}})H$. Hence $v_s$ is a partial isometry 
with initial space $\pi (E_{s^{*}})H$, hence final space $\pi (E_
s)H$ 
(since $v_s^{*}=\pi (w_{s^{*},s}^{*})v_{s^{*}}$).

The covariance condition is satisfied since if $a\in E_{s^{*}}$ then 
$$\eqalign{v_s\pi (a)v_s^{*}&=v_s\pi (a)\pi (w_{s^{*},s}^{*})v_{s^{
*}}\cr
&=\slim_{\mu ,\lambda}\Pi (u_{\lambda}\delta_s*aw_{s^{*},s}^{*}\delta_
e*\beta_s^{-1}(u_{\mu})\delta_{s^{*}})\cr
&=\slim_{\mu ,\lambda}\Pi (\beta_s(\beta_s^{-1}(u_{\lambda})aw_{s^{
*},s}^{*})\delta_s*\beta_s^{-1}(u_{\mu})\delta_{s^{*}})\cr
&=\slim_{\mu ,\lambda}\Pi (u_{\lambda}\beta_s(aw_{s^{*},s}^{*})u_{
\mu}w_{s,s^{*}}\delta_{ss^{*}})\cr
&=\Pi (\beta_s(a)\delta_e)\,,\cr}
$$
since $\{\beta_s^{-1}(u_{\mu})\}$ is an approximate identity for $
E_{s^{*}}$. 
\eop
 
\prop The correspondence $(\pi ,v,H)\leftrightarrow (\pi\times v,
H)$ is a bijection 
between covariant representations of $(A,S,\beta ,w)$ and 
nondegenerate representations of $A\times_{\beta ,w}S$.  

\label {hivz} 

\proof Let $(\tilde{\pi },\tilde {v},H)$ be a covariant representation of 
$(A,S,\beta ,w)$, and let $(\pi ,v)$ be the covariant representation 
associated to $\tilde{\pi}\times\tilde {v}$ by Proposition \cite{slimref}.  Then 
$$\eqalign{\pi (a)&=\big(\tilde{\pi}\times\tilde {v}\big)(a\delta_
e)=\tilde{\pi }(a)\tilde {v}_e=\tilde{\pi }(a)\cr
v_s&=\slim_{\lambda}\big(\tilde{\pi}\times\tilde {v}\big)(u_{\lambda}
\delta_s)=\slim_{\lambda}\tilde{\pi }(u_{\lambda})\tilde {v}_s=\tilde {
v}_s\,,\cr}
$$
since $\tilde{\pi }(u_{\lambda})$ converges 
strongly to the projection of $H$ onto $\tilde{\pi }(E_s)H$.  On the other hand, if 
$\Pi$ is a representation of $A\times_{\beta}S$ and $(\pi ,v)$ is the covariant 
representation associated to $\Pi$, then for $a\in E_s$ we have 
$$\eqalign{\big(\pi\times v\big)(a\delta_s)&=\pi (a)v_s=\Pi (a\delta_
e)\slim_{\lambda}\Pi (u_{\lambda}\delta_s)\cr
&=\slim_{\lambda}\Pi (a\delta_e*u_{\lambda}\delta_s)\cr
&=\slim_{\lambda}\Pi (au_{\lambda}\delta_s)=\Pi (a\delta_s).\cr}
$$
Thus the correspondence is a bijection.
\eop

\newsection{Connection with twisted partial actions}

\uj We have seen in [Sie] and [Ex3] that there is a close 
connection between partial actions of [McC] and 
untwisted inverse semigroup actions.  There is a similar 
connection between {\it twisted partial actions\/} of [Ex2] and 
Busby-Smith twisted inverse semigroup actions, which is 
the topic of this section. Recall the definition of a 
twisted partial action from [Ex2].

\definition A {\it twisted partial action\/} of a group $G$ on a 
$C^{*}$-algebra $A$ is a pair $(\alpha ,u)$, where for all $s\in 
G$, 
$\alpha_s:D_{s^{-1}}\to D_s$ is a partial automorphism of $A$, and for all 
$r$, $s\in G$, $u_{r,s}$ is a unitary multiplier of $D_rD_{rs}$, such that 
for all $r$, $s$, $t\in G$ we have 
\item{(a)}{$D_e=A$, and $\alpha_e$ is the identity automorphism 
of $A$;} 
\item{(b)}{$\alpha_r(D_{r^{-1}}D_s)=D_rD_{rs}$;} 
\item{(c)}{$\alpha_r(\alpha_s(a))=u_{r,s}\alpha_{rs}(a)u_{r,s}^{*}$ for all 
$a\in D_{s^{-1}}D_{s^{-1}r^{-1}}$;} 
\item{(d)}{$u_{e,t}=u_{t,e}=1_{M(A)}$;} 
\item{(e)}{$\alpha_r(au_{s,t})u_{r,st}=\alpha_r(a)u_{r,s}u_{rs,t}$ for all 
$a\in D_{r^{-1}}D_sD_{st}$;} 
\label{partialaction} 

\uj Recall from [Sie], that if $\alpha$ is a partial action of $G$ on 
$A$ then the partial automorphism $\alpha_{s_1}\cdots\alpha_{s_n}$ has domain 
$D_{s_n^{-1}}D_{s_n^{-1}s_{n-1}^{-1}}\cdots D_{s_n^{-1}\cdots s_1^{
-1}}$ and range $D_{s_1}D_{s_1s_2}\cdots D_{s_1\cdots s_n}$ for 
all $s_1,\ldots ,s_n\in G$.  A similar proof shows that this is also 
true for twisted partial actions, even if $\alpha_{s_i}$ is 
replaced by $\alpha_{s_i^{-1}}^{-1}$ for some $i$.  

Recall from [Ex3] that for a group $G$, the {\it associated }
{\it inverse semigroup} $S(G)$ has elements written in  
canonical form $[g_1][g_1^{-1}]\cdots [g_m][g_m^{-1}][s]$, where 
$g_1,\dots,g_n,s\in G$, and the order of the $[g_i][g_i^{-1}]$ terms is 
irrelevant. Multiplication and inverses are defined by
$$\eqalign{&[g_1][g_1^{-1}]\cdots [g_m][g_m^{-1}][s]\cdot [h_1][h_
1^{-1}]\cdots [h_m][h_m^{-1}][t]\cr
&\quad =[g_1][g_1^{-1}]\cdots [g_m][g_m^{-1}][s][s^{-1}][sh_1][(s
h_1)^{-1}]\cdots [sh_m][(sh_m)^{-1}][st]\cr}
$$
and
$$([g_1][g_1^{-1}]\cdots [g_m][g_m^{-1}][s])^{*}=[s^{-1}g_m][(s^{
-1}g_m)^{-1}]\cdots [s^{-1}g_1][(s^{-1}g_1)^{-1}][s^{-1}]\,.$$
The idempotents of $S(G)$ are the elements in the form 
$[g_1][g_1^{-1}]\cdots [g_m][g_m^{-1}][e]$, where $e$ is the identity of $
G$. The 
next theorem shows that every twisted partial action of $G$ 
determines a Busby-Smith twisted action of $S(G)$.   

\mark{parcon}
\theorem Let $(A,G,\alpha ,u)$ be a twisted partial action.  For all 
$p=[g_1][g_1^{-1}]\cdots [g_m][g_m^{-1}][s]$ and $q=[h_1][h_1^{-1}
]\cdots [h_n][h_n^{-1}][t]$ 
in $S(G)$, define $E_p=D_{g_1}\cdots D_{g_m}D_s$ and let
$$\beta_p=\alpha_{g_1}\alpha_{g_1}^{-1}\cdots\alpha_{g_m}\alpha_{
g_m}^{-1}\alpha_s\,.$$
Also let 
$$w_{p,q}=1_{M(E_{pq})}u_{s,t}\,.$$
Then $(A,S(G),\beta ,w)$ is a Busby-Smith twisted action.

\proof Throughout the proof, let $p$ and $q$ be in the form 
$p=[g_1][g_1^{-1}]\cdots [g_m][g_m^{-1}][s]$ and $q=[h_1][h_1^{-1}
]\cdots [h_n][h_n^{-1}][t]$. 
First note that for all $p\in S(G)$, $\beta_p$ is an isomorphism 
between the closed ideals $E_{p^{*}}$ and $E_p$ of $A$.
Also note 
that $w_{p,q}$ is a unitary multiplier of $E_{p,q}$, since $u_{s,
t}$ is a 
unitary multiplier of $D_sD_{st}$ and 
$E_{pq}=D_{g_1}\cdots D_{g_m}D_sD_{sh_1}\cdots D_{sh_n}D_{st}$ is contained in $
D_sD_{st}$.  
If $e$ is the identity of $G$ then $E_{[e]}=D_e=A$, verifying 
Definition \cite{taction}(a).
For all 
$p,q\in S(G)$ we have 
$$\beta_p\beta_q=\alpha_{g_1}\alpha_{g_1}^{-1}\cdots\alpha_{g_m}\alpha_{
g_m}^{-1}\alpha_s\alpha_{h_1}\alpha_{h_1}^{-1}\cdots\alpha_{h_n}\alpha_{
h_n}^{-1}\alpha_t\,\,,$$
which is the restriction of $\alpha_s\alpha_t$ to 
$$D_{t^{-1}}D_{t^{-1}h_n}\cdots D_{t^{-1}h_1}D_{t^{-1}s^{-1}}D_{t^{
-1}s^{-1}g_m}\cdots D_{t^{-1}s^{-1}g_1}=E_{(pq)^{*}}\,.$$
On the other hand 
$$\beta_{pq}=\alpha_{g_1}\alpha_{g_1}^{-1}\cdots\alpha_{g_m}\alpha_{
g_m}^{-1}\alpha_s\alpha_s^{-1}\alpha_{sh_1}\alpha_{sh_1}^{-1}\cdots
\alpha_{sh_n}\alpha_{sh_n}^{-1}\alpha_{st}\,,$$
which is the restriction of $\alpha_{st}$ to $E_{(pq)^{*}}$.  Thus, 
Definition \cite{taction}(b) follows from 
Definition \cite{partialaction}(c), since $E_{(pq)^{*}}\subset D_{
t^{-1}}D_{t^{-1}s^{-1}}$ .  To 
check Definition \cite{taction}(c), note that if 
$f=[g_1][g_1^{-1}]\cdots [g_m][g_m^{-1}][e]$ is an idempotent in S(G), then 
$w_{f,q}=1_{M(E_{fq})}u_{e,t}=1_{M(E_{fq})}$ for all $q\in S(G)$.  Similarly, 
$w_{q,f}=1_{M(E_{qf})}$ for all $q\in S(G)$.  Finally to check 
Definition
\cite{taction}(d), let $k=[f_1][f_1^{-1}]\cdots [f_l][f_l^{-1}][r
]$, $p$ and $q$ be 
arbitrary elements of $S(G),$ and fix 
$$a\in E_{k^{*}}E_{pq}=D_{r^{-1}}D_{r^{-1}f_l}\cdots D_{r^{-1}f_1}
D_{g_1}\cdots D_{g_n}D_sD_{sh_1}\cdots D_{sh_n}D_{st}\,.$$
Then $a\in D_{r^{-1}}D_sD_{st}$, and so we have   
$$\eqalign{\beta_k(aw_{p,q})w_{k,pq}&=\beta_k(a1_{M(E_{pq})}u_{s,
t})1_{M(E_{kpq})}u_{r,st}\cr
&=\alpha_r(au_{s,t})u_{r,st}=\alpha_r(a)u_{r,s}u_{rs,t}\cr
&=\alpha_r(a)1_{M(E_{kp})}u_{r,s}1_{M(E_{kpq})}u_{rs,t}\cr
&=\beta_k(a)w_{k,p}w_{kp,q}\,.\cr}
$$
\Eop

\uj This process works the other way too. Starting 
with a Busby-Smith twisted action $\beta$ of $S(G)$, the 
restriction of $\beta$ to the canonical image of $G$ in $S(G)$ 
gives a twisted partial action:

\theorem Let $(A,S(G),\beta ,w)$ be a Busby-Smith twisted 
action. If $\alpha_s=\beta_{[s]}$ and $u_{s,t}=w_{[s],[t]}$ for all 
$s$,$t\in G$, then $(A,G,\alpha ,u)$ is a twisted partial action.

\proof Let $D_s$ be the range of $\alpha_s$ for all $s\in G$. If $
e$ is the 
identity of $G$, then $D_e=E_{[e]}=A$. We have 
$$\eqalign{\alpha_r(D_{r^{-1}}D_s)&=\hbox{\rm im}\,\beta_{[r]}\beta_{
[s]}=\hbox{\rm im}\,\beta_{[r][s]}=\hbox{\rm im}\,\beta_{[r][r^{-
1}][rs]}\cr
&=\hbox{\rm im}\,\beta_{[r]}\beta_{[r^{-1}]}\beta_{[rs]}=E_{[r]}E_{
[rs]}=D_rD_{rs}\,,\cr}
$$
which verifies Definition \cite{partialaction}(b). Similar 
calculations can be used to verify the other conditions in 
Definition \cite{partialaction}.
\eop

\uj It is clear that the processes above are the 
inverses of each other. Thus we have a bijective 
correspondence between twisted partial actions of $G$ and 
Busby-Smith twisted actions of $S(G)$. 
Although the theory of crossed products by twisted 
partial actions has not been developed, we believe that the 
crossed products of corresponding twisted partial 
actions mentioned in [Ex2] and Busby-Smith twisted actions are isomorphic. 
This is, in fact, the case for untwisted actions, since 
the semigroup action constructed in Theorem \cite{parcon} 
factors through the semigroup action giving the 
isomorphic crossed product constructed in [Sie]. We plan to pursue this 
in an upcoming paper. 

It would be 
interesting to know if there is a similar correspondence 
for the Green twisted actions defined in Section 6 below.
It would first be necessary to find a satisfactory 
definition of Green twisted partial actions of a group.

\newsection{Exterior equivalence}

\uj Exterior equivalence is defined by Packer and Raeburn 
for twisted group actions in [PR].  We extend this 
definition to Busby-Smith twisted inverse semigroup actions.  

\definition Two Busby-Smith twisted actions $(\alpha ,u)$ 
and $(\beta ,w)$ of $S$ on $A$ are {\it exterior equivalent\/} if for all 
$s\in S$ there is a unitary multiplier $V_s$ of $E_s$ such that 
for all $s,t\in S$ 
\item{(a)}{$\beta_s=\hbox{\rm Ad }V_s\circ\alpha_s$;} 
\item{(b)}{$w_{s,t}=V_s\alpha_s(1_{M(E_{s^{*}})}V_t)u_{s,t}V_{st}^{
*}$.} 

\uj We say that the exterior equivalence is {\it implemented }
by $V$.  Notice that for all $s\in S$, $\alpha_s$ and $\beta_s$ have to have 
the same domain and range.  Also note that 
$1_{M(E_{s^{*}})}V_t\in E_{s^{*}}^{**}$ since $E_{s^{*}}^{**}=1_{
M(E_{s^{*}})}A^{**}$.  Hence we can 
evaluate $\alpha_s(1_{M(E_{s^{*}})}V_t)$ if we extend $\alpha_s$ to the double 
dual of $E_{s^{*}}$.  

\prop Exterior equivalence is an equivalence relation.

\proof Let $V$ implement an exterior equivalence between 
the Busby-Smith twisted actions $(\alpha ,u)$ and $(\beta ,w)$.  Then taking 
adjoints we have $\alpha_s=\hbox{\rm Ad }V_s^{*}\circ\beta_s$ and 
$u_{s,t}=\alpha_s(1_{M(E_{s^{*}})}V_t^{*})V_s^{*}w_{s,t}V_{st}=V_
s^{*}\beta_s(1_{M(E_{s^{*}})}V_t^{*})w_{s,t}V_{st}$, so 
$V^{*}$ implements an exterior equivalence between $(\beta ,w)$ and 
$(\alpha ,u)$, showing symmetry.  Reflexivity is clear by letting 
$V_s=1_{M(E_s)}$ for all $s\in S$.  To show transitivity let $V$ 
implement an exterior equivalence between $(\alpha ,u)$ and 
$(\beta ,w)$, and $X$ implement an exterior equivalence between 
$(\beta ,w)$ and $(\gamma ,z)$.  Then 
$$\eqalign{\gamma_s&=\hbox{\rm Ad }X_sV_s\circ\alpha_s\,,\cr
z_{s,t}&=X_s\beta_s(1_{M(E_{s^{*}})}X_t)w_{s,t}X_{st}^{*}\cr
&=X_sV_s\alpha_s(1_{M(E_{s^{*}})}X_t)V_s^{*}V_s\alpha_s(1_{M(E_{s^{
*}})}V_t)u_{s,t}V_{st}^{*}X_{st}^{*}\cr
&=X_sV_s\alpha_s(1_{M(E_{s^{*}})}X_tV_t)u_{s,t}(X_{st}V_{st})^{*}\,
,\cr}
$$
which shows that $XV$ implements an exterior 
equivalence between $(\alpha ,u)$ and $(\gamma ,z)$.
\hbox{\hskip .1in}\eop

\prop If $(\alpha ,u)$ and $(\beta ,w)$ are exterior equivalent 
Busby-Smith twisted inverse semigroup actions of $S$ on 
$A$, then $A\times_{\alpha ,u}S$ and $A\times_{\beta ,w}S$ are isomorphic.  

\proof Suppose $V$ implements an exterior equivalence 
between $(\alpha ,u)$ and $(\beta ,w).$ Let $(\pi ,z)$ be a covariant 
representation of $(\beta ,w)$, and define $v_s=\pi (V_s^{*})z_s$.  We 
show that $(\pi ,v)$ is a covariant representation of $(\alpha ,u
)$.  
If $s\in S$ and $a\in E_{s^{*}}$ then 
$$\pi (\alpha_s(a))=\pi (V_s^{*}\beta_s(a)V_s)=\pi (V_s^{*})z_s\pi 
(a)z_s^{*}\pi (V_s)=v_s\pi (a)v_s^{*}.$$
We also have 
$$\eqalign{v_sv_t&=\pi (V_s^{*})z_s\pi (1_{M(E_{s^{*}})}V_t^{*})z_
t\cr
&=\pi (V_s^{*})\pi (\beta_s(1_{M(E_{s^{*}})}V_t^{*}))z_sz_t\cr
&=\pi (V_s^{*})\pi (V_s)\pi (\alpha_s(1_{M(E_{s^{*}})}V_t^{*}))\pi 
(V_s^{*})\pi (w_{s,t})z_{st}\cr
&=\pi (u_{s,t})\pi (V_{st}^{*})z_{st}=\pi (u_{s,t})v_{st}.\cr}
$$
The other conditions of Definition \cite{covrepdef} are 
clearly satisfied.  Note that the images of $\pi\times v$ and $\pi
\times z$ 
are the same, that is, $\big(\pi\times v\big)(A\times_{\alpha ,u}
S)=\big(\pi\times z\big)(A\times_{\beta ,w}S)$.  

Now suppose that $\pi\times z$ is a faithful representation of 
$A\times_{\beta ,w}S$.  Then $\Phi =$$(\pi\times z)^{-1}\circ (\pi
\times v):A\times_{\alpha ,u}S\to A\times_{\beta ,w}S$ is a 
surjective homomorphism.  Similarly, starting with a 
faithful representation $\rho\times l$ of $A\times_{\alpha ,u}S$, we can find a 
covariant representation $(\rho ,m)$ of $(\beta ,w)$ such that 
$\Psi =(\rho\times l)^{-1}\circ (\rho\times m):A\times_{\beta ,w}
S\to A\times_{\alpha ,u}S$ is a 
homomorphism.  We are going to show that $\Psi\circ\Phi$ is the 
identity map, which implies that $\Phi$ is an isomorphism.  
For $a\in E_s$ we have 
$$\eqalign{\Psi\circ\Phi (a\delta_s)&=\Psi\circ (\pi\times z)^{-1}
(\pi (a)v_s)=\Psi\circ (\pi\times z)^{-1}(\pi (a)\pi (V_s^{*})z_s
)\cr
&=\Psi\circ (\pi\times z)^{-1}(\pi (aV_s^{*})z_s)=\Psi (aV_s^{*}\delta_
s)\cr
&=(\rho\times l)^{-1}(\rho (aV_s^{*})m_s)=(\rho\times l)^{-1}(\rho 
(aV_s^{*})\rho (V_s)l_s)\cr
&=(\rho\times l)^{-1}(\rho (aV_s^{*}V_s)l_s)=a\delta_s\,.\cr}
$$
\Eop

\example If we have two order-preserving 
cross-sections $c,d:S\to T$ in Proposition \cite{tactionexample}, 
then the corresponding Busby-Smith twisted actions $(\beta ,w)$ 
and $(\gamma ,z)$ defined by $c$ and $d$ respectively, 
are exterior equivalent, and so $C^{*}(N)\times_{\beta ,w}S$ and 
$C^{*}(N)\times_{\gamma ,z}S$ are isomorphic.  

To see this define $V_s=d(s)c(s)^{*}$ for all $s\in S$. Then we have
$$\eqalign{\gamma_s=\hbox{\rm Ad }d(s)&=\hbox{\rm Ad }d(s)\circ\hbox{\rm Ad }
c(s)^{*}\circ\hbox{\rm Ad }c(s)\cr
&=\hbox{\rm Ad }d(s)c(s)^{*}\circ\hbox{\rm Ad }c(s)=\hbox{\rm Ad }
V_s\circ\beta_{s\,,}\cr}
$$
and
$$\eqalign{&V_s\beta_s(1_{M(E_{s^{*}})}V_t)w_{s,t}V_{st}^{*}\cr
&\quad =d(s)c(s)^{*}c(s)1_{M(E_{s^{*}})}d(t)c(t)^{*}c(s)^{*}c(s)c
(t)c(st)^{*}(d(st)c(st)^{*})^{*}\cr
&\quad =d(s)d(t)d(st)^{*}=z_{s,t}\cr}
$$
showing that $V$ implements an exterior equivalence 
between $(\beta ,w)$ and $(\gamma ,z)$.   \eop

\newsection{Green twisted actions}

\uj Green studies another type of twisted group action in 
[Gre], and here we adapt this to inverse semigroups:  

\definition Let $A$ be a $C^{*}$-algebra, let $S$ be a unital inverse 
semigroup with idempotent semilattice $E$, and let $N$ be a 
normal Clifford subsemigroup of $S$.  A {\it Green twisted action }
of $(S,N)$ on $A$ is a pair $(\gamma ,\tau )$, where $\gamma$ is an inverse 
semigroup action of $S$ on $A$, that is, a semigroup 
homomorphism $s\mapsto (\gamma_s,E_{s^{*}},E_s):S\to\hbox{\rm PAut }
(A)$ with $E_e=A$, 
and for all $n\in N$, $\tau_n$ is a unitary multiplier of $E_n$, such 
that for all $n,l\in N$ we have 
\item{(a)}{$\gamma_n=\hbox{\rm Ad }\tau_n$;} 
\item{(b)}{$\gamma_s(\tau_n)=\tau_{sns^{*}}$ for all $s\in S$ with $
n^{*}n\leq s^{*}s$;} 
\item{(c)}{$\tau_n\tau_l=\tau_{nl}$.} 

\noindent We call $\tau$ the twisting map, and we also refer to 
$(A,S,N,\gamma ,\tau )$ as a Green twisted action.  

\label{gtactiondef}

\example Let $S$, $E$, $N$ be as in Definition 
\cite{gtactiondef}.  Define 
$E_s=\hbox{\rm $\overline {\hbox{\rm span}}$}\,\cup \{[f]:f\leq ss^{
*}\hbox{\rm \ for some idempotent }f\in E\}\subset C^{*}(N)$.  
For all $s\in S$ define $\gamma_s:E_{s^{*}}\to E_s$ by $\gamma_s=\hbox{\rm Ad }
s$, and for all 
$n\in N$ define $\tau_n=n$.  Then $(C^{*}(N),S,N,\gamma ,\tau )$ is a Green 
twisted action.  

\label{gtactionexample}

\definition Let $(A,S,N,\gamma ,\tau )$ be a Green twisted action.  A 
{\it covariant representation\/} of the Green twisted action 
$(\gamma ,\tau )$ is a covariant representation $(\pi ,u)$ of the action $
\gamma$ 
that {\it preserves the twist} $\tau$, that is, $\pi$ is a nondegenerate 
representation of $A$ on the Hilbert space $H$, and 
$u:S\to B(H)$ is a multiplicative map such that 
\item{(a)}{$u_s\pi (a)u_s^{*}=\pi (\gamma_s(a))$ for all $a\in E_{
s^{*}}$;} 
\item{(b)}{$u_s$ is a partial isometry with initial space 
$\pi (E_{s^{*}})H$ and final space $\pi (E_s)H$;} 
\item{(c)}{$u_n=\pi (\tau_n)$ for all $n\in N$.} 

\noindent
To evaluate $\pi (\tau_n)$ we again extend $\pi$ to the enveloping von 
Neumann algebra $A^{**}$ of $A$.

\label{gcovrepdef}

\definition Let $(A,S,N,\gamma ,\tau )$ be a Green twisted action.  
The {\it Green twisted crossed product} $A\times_{\gamma ,\tau}S$ is the 
quotient of $A\times_{\gamma}S$ by the ideal 
$$I_{\tau}=\cap \{\hbox{\rm ker}\,(\pi\times u):(\pi ,u)\hbox{\rm \ is a covariant representation of }
(\gamma ,\tau )\}.$$
  
\uj Our definitions of Green twisted inverse semigroup 
action and Green twisted crossed product are 
generalizations of inverse semigroup action and crossed 
product defined in [Sie].  Every action $\beta$ of $S$ may be 
regarded (trivially) as a Green twisted inverse semigroup 
action by taking $\tau_n=1_{M(E_n)}$ for all $n$ in the idempotent 
semilattice $E$ of $S$, and we have $A\times_{\beta}S\cong A\times_{
\beta ,\tau}S$.  

\theorem There is a bijective correspondence between 
nondegenerate representations of the Green twisted 
crossed product and covariant representations of the 
Green twisted action.  

\proof Let $\Psi :A\times_{\gamma}S\to A\times_{\gamma ,\tau}S$ be the quotient map.  If $
\Pi$ is 
a representation of $A\times_{\gamma ,\tau}S$ then $\Pi\circ\Psi$ is a 
representation of $A\times_{\gamma}S$ and so $\Pi\circ\Psi =\pi\times 
u$ for some 
covariant representation $(\pi ,u)$ of $\gamma$ by Proposition \cite{hivz}.  
We show that $(\pi ,u)$ 
preserves the twist.  If $a\in E_n$ then for every covariant 
representation $(\rho ,z)$ of $(\gamma ,\tau )$ we have 
$\big(\rho\times z\big)(a\tau_n\delta_e-a\delta_n)=\rho (a)(\rho 
(\tau_n)-z_n)=0$ and so 
$$\pi (a)(\pi (\tau_n)-u_n)=\big(\pi\times u\big)(a\tau_n\delta_e
-a\delta_n)=\Pi\circ\Psi (a\tau_n\delta_e-a\delta_n)=\Pi (0)=0.$$
Thus $\big(\pi (\tau_n)^{*}-u_n^{*}\big)\pi (E_n)H=\{0\}$. Since $
\pi (\tau_n)^{*}$ and $u_n^{*}$ are 
partial isometries with initial space $\pi (E_n)H$, we must 
have $\pi (\tau_n)-u_n=0$ and so $(\pi ,u)$ is a covariant 
representation of $(\gamma ,\tau )$.

On the other hand, if $(\pi ,u)$ is a covariant representation 
of the twisted action $(\gamma ,\tau )$, then $I_{\tau}\subset\hbox{\rm ker }
(\pi\times u)$, and so 
there is a unique representation $\pi\times_{\tau}u$ of $A\times_{
\gamma ,\tau}S$ defined 
by $\pi\times_{\tau}u(\Psi (x))=\big(\pi\times u\big)(x)$.   

For the uniqueness, note that Proposition \cite{hivz}, in the absence 
of the Busby-Smith twist, gives a bijection between nondegenerate 
representations $\Pi$ of $A\times\gamma S$ and covariant representations $
(\pi ,u)$ of 
$(A,S,\gamma )$. By the above argument, $\Pi$ kills the ideal $I_{
\tau}$ if and only if 
$(\pi ,u)$ preserves the twist $\tau$, so we are done. 
\eop

\newsection{Connection between Busby-Smith and Green 
twisted actions} 

\uj There is a close connection between Busby-Smith 
twisted and Green twisted actions just like in the group 
case [PR].  In this section we show that starting with a
Busby-Smith twisted action we can construct a Green 
twisted action with the same crossed product.  
Conversely starting with a Green twisted action we 
construct a Busby-Smith twisted action with the same 
crossed product.  This latter construction depends on the 
existence of an order-preserving cross-section, 
suggesting that Green twisted actions are ``more 
general'' than Busby-Smith twisted actions.  This may 
seem surprising since in the discrete group case there is 
an essentially bijective correspondence between 
Busby-Smith twisted and Green twisted actions.  

\theorem Let $(A,S,N,\gamma ,\tau )$ be a Green twisted  
action, and suppose there is an order-preserving cross-section 
$c:S/N\to S$.  For $q,r\in S/N$ define 
$$\beta_q=\gamma_{c(q)}\quad\hbox{\rm and}\quad w_{q,r}=\tau_{c(q
)c(r)c(qr)^{*}}\,.$$
Then $(A,S/N,\beta ,w)$ is a Busby-Smith twisted action, and 
the crossed products $A\times_{\gamma ,\tau}S$ and $A\times_{\beta 
,w}S/N$ are 
isomorphic. 
\label{connectionthm1} 

\proof First we show that $(\beta ,w)$ is a Busby-Smith twisted 
action.  If $q$, $r\in S/N$ then $w_{q,r}$ is a 
unitary multiplier of $E_{c(q)c(r)}=E_{c(q)c(r)c(r)^{*}c(q)^{*}}$ since 
$\tau_{c(q)c(r)c(qr)^{*}}$ is a unitary multiplier of 
$E_{c(q)c(r)c(qr)^{*}}=E_{c(q)c(r)c(qr)^{*}c(qr)c(r)^{*}c(q)^{*}}
=E$$_{c(q)c(r)c(r)^{*}c(q)^{*}}$.  
It is clear that $E_{[e]}=E_e=A$. For $q$, $r\in S/N$ we have 
$$\eqalign{\beta_q\beta_r&=\gamma_{c(q)}\gamma_{c(r)}=\gamma_{c(q
)c(r)}=\gamma_{c(q)c(r)c(qr)^{*}c(qr)}\cr
&=\hbox{\rm Ad }\tau_{c(q)c(r)c(qr)^{*}}\circ\gamma_{c(qr)}=\hbox{\rm Ad }
w_{q,r}\circ\beta_{qr}\,,\cr}
$$
which verifies Definition \cite{taction}(b).  To check 
Definition \cite{taction}(c) let $f\in E$ and $q\in S/N$.  Then 
we have 
$$\eqalign{w_{[f],q}&=\tau_{fc(q)c([f]q)^{*}}=1_{M(E_{[f]qq^{*}})}
=1_{M(E_{[f]q})}\cr}
.$$
To see this notice that since $c$ is order-preserving and 
$q\ge [f]q$, we have $c(q)\ge c([f]q)$. Hence 
$fc(q)c([f]q)^{*}=fc([f]q)c([f]q)^{*}$ is an idempotent in $[f]qq^{
*}$. 
Similarly $w_{q,[f]}=1_{M(E_{q[f]})}.$ It remains to check 
Definition
\cite{taction}(d). If $p$, $q$, $r\in S/N$ and $a\in E_{c(p^{*})}
E_{c(qr)}$, 
then
$$\eqalign{\beta_p(aw_{q,r})w_{p,qr}&=\gamma_{c(p)}(a\tau_{c(q)c(
r)c(qr)^{*}})\tau_{c(p)c(qr)c(pqr)^{*}}\cr
&=\gamma_{c(p)}(a)\tau_{c(p)c(q)c(r)c(qr)^{*}c(p)^{*}c(p)c(qr)c(p
qr)^{*}}\cr
&=\gamma_{c(p)}(a)\tau_{c(p)c(q)c(r)c(pqr)^{*}}\cr
&=\gamma_{c(p)}(a)\tau_{c(p)c(q)c(pq)^{*}c(pq)c(r)c(pqr)^{*}}\cr
&=\beta_p(a)w_{p,q}w_{pq,r}\,.\cr}
$$
 
Next we investigate the connection between $(\gamma ,\tau )$ and 
$(\beta ,w)$.  Let $(\pi ,v)$ be a covariant representation of $(
\beta ,w)$.  
Define $u_s=\pi (\tau_{sc([s])^{*}})v_{[s]}$.  We show that $(\pi 
,u)$ is a 
covariant representation of $(\gamma ,\tau )$.  First notice that $
u$ is 
a homomorphism since for $s$, $t\in S$ we have 
$$\eqalign{u_su_t&=\pi (\tau_{sc([s])^{*}})v_{[s]}\pi (\tau_{tc([
t])^{*}})v_{[t]}\cr
&=\pi (\tau_{sc([s])^{*}})\pi (\beta_{[s]}(\tau_{tc([t])^{*}}))v_{
[t]}v_{[t]}\cr
&=\pi (\tau_{sc([s])^{*}}\tau_{c([s])tc([t])^{*}c([s])^{*}})\pi (
w_{[s],[t]})v_{[st]}\cr
&=\pi (\tau_{stc([t])^{*}c([s])^{*}})\pi (\tau_{c([s])c([t])c([st
])^{*}})v_{[st]}\cr
&=\pi (\tau_{stc([st])^{*}})v_{[st]}=u_{st}\,.\cr}
$$
To check the covariance condition let $s\in S$ and $a\in E_{s^{*}}$, 
then
$$\eqalign{u_s\pi (a)u_s^{*}&=\pi (\tau_{sc([s])^{*}})v_{[s]}\pi 
(a)v_{[s]}^{*}\pi (\tau_{sc([s])^{*}}^{*})\cr
&=\pi (\gamma_{sc([s])^{*}}(\beta_{[s]}(a)))\quad\hbox{\rm (Definition \cite{gtactiondef}(a) )}\cr
&=\pi (\gamma_{sc([s])^{*}c([s])}(a))=\pi (\gamma_s(a))\cr}
$$
where we used the fact that $sc([s])^{*}\in N$.  It is clear 
that $u_s$ has the right initial and final spaces.  The 
following calculation shows that $(\pi ,u)$ preserves the 
twist.  
$$\eqalign{u_n&=\pi (\tau_{nc([n])^{*}})v_{[n]}=\pi (\tau_n)v_{[n
]}\cr
&=\pi (\tau_n)P_{\pi (E_{[n]})H}\quad\hbox{\rm (Proposition \cite{vproperty}(a) )}\cr
&=\pi (\tau_n)P_{\pi (E_n)H}=\pi (\tau_n)\,.\cr}
$$
We show that $\hbox{\rm im}\,(\pi\times_{\tau}u)=\hbox{\rm im}\,(
\pi\times v)$. For $a\in E_s$ we have  
$$\pi (a)u_s=\pi (a)\pi (\tau_{sc([s])^{*}})v_{[s]}=\pi (a\tau_{s
c([s])^{*}})v_{[s]}\in\hbox{\rm im}\,(\pi\times v)$$
and so $\hbox{\rm im}\,(\pi\times_{\tau}u)\subset\hbox{\rm im}\,(
\pi\times v)$. On the other hand,
$$\pi (a)v_{[s]}=\pi (a\tau_{sc([s])^{*}}^{*})\pi (\tau_{sc([s])^{
*}})v_{[s]}=\pi (a\tau_{sc([s])^{*}}^{*})u_s\in\hbox{\rm im}\,(\pi
\times_{\tau}u)$$
and so $\hbox{\rm im}\,(\pi\times_{\tau}u)\supset\hbox{\rm im}\,(
\pi\times v)$.

Next let $(\pi ,u)$ be a covariant representation of $(\gamma ,\tau 
)$. 
Define $v_q=u_{c(q)}$. We show that $(\pi ,v)$ is a covariant 
representation of $(\beta ,w)$. If $a\in E_{q^{*}}$ then
$$\eqalign{\pi (\beta_q(a))=\pi (\gamma_{c(q)}(a))=u_{c(q)}\pi (a
)u_{c(q)}^{*}=v_q\pi (a)v_{q^{*}}\cr}
\,,$$
which shows that the covariance condition holds. We also 
have 
$$\eqalign{v_qv_r&=u_{c(q)}u_{c(r)}=u_{c(q)c(r)c(qr)^{*}}u_{c(qr)}\cr
&=\pi (\tau_{c(q)c(r)c(qr)^{*}})u_{c(qr)}=\pi (w_{q,r})v_{qr}\,.\cr}
$$
It is clear that $v_q$ has the right initial and final spaces.  
It is also clear that $\hbox{\rm im}\,(\pi\times_{\tau}u)\supset\hbox{\rm im}\,
(\pi\times v)$ since for $a\in E_{[s]}$ we 
have $\pi (a)v_{[s]}=\pi (a)u_{c([s])}\in\hbox{\rm im}\,(\pi\times_{
\tau}u)$. 

Now let $\pi\times v$ be a faithful representation of $A\times_{\beta 
,w}S/N$.  
Then by the above correspondence, $\pi\times_{\tau}u$ is a 
representation of $A\times_{\gamma ,\tau}S$ such that 
$\hbox{\rm im}\,(\pi\times v)=\hbox{\rm im}\,(\pi\times_{\tau}u)$, and 
so we have a surjective homomorphism 
$\Phi =(\pi\times v)^{-1}\circ (\pi\times_{\tau}u):A\times_{\gamma 
,\tau}S\to A\times_{\beta ,w}S/N$.  Similarly, 
starting with a faithful representation $\rho\times_{\tau}z$ of $
A\times_{\gamma ,\tau}S$ 
the representation $\rho\times (z\circ c)$ can be used to find a 
homomorphism 
$\Psi =(\rho\times_{\tau}z)^{-1}\circ (\rho\times (z\circ c)):A\times_{
\beta ,w}S/N\to A\times_{\gamma ,\tau}S$. We show 
that $\Psi\circ\Phi$ is the identity map, hence $\Phi$ is an 
isomorphism. For $a\in E_s$ we have  
$$\eqalign{\Psi\circ\Phi (a\delta_s)&=\Psi\circ (\pi\times v)^{-1}
(\pi (a)u_s)\cr
&=\Psi\circ (\pi\times v)^{-1}(\pi (a)\pi (\tau_{sc([s])^{*}})v_{
[s]})\cr
&=\Psi (a\tau_{sc([s])^{*}}\delta_{[s]})\cr
&=(\rho\times_{\tau}z)^{-1}\circ\big(\rho\times (z\circ c)\big)(a
\tau_{sc([s])^{*}}\delta_{[s]})\cr
&=(\rho\times_{\tau}z)^{-1}(\rho (a)\rho (\tau_{sc([s])^{*}})z_{c
[s]})\cr
&=(\rho\times_{\tau}z)^{-1}(\rho (a)z_{sc([s])^{*}}z_{c[s]})\cr
&=(\rho\times_{\tau}z)^{-1}(\rho (a)z_s)=a\delta_s\,.\cr}
$$
\Eop

\prop If we use a different order-preserving cross-section 
$b:S/N\to S$ in Theorem \cite{connectionthm1} then we get 
an exterior equivalent Busby-Smith twisted action $(\alpha ,u)$. 

\proof The exterior equivalence is implemented by 
$V_q=\tau_{c(q)b(q)^{*}}$ because for $q\in S/N$ we have 
$$\beta_q=\gamma_{c(q)}=\gamma_{c(q)b(q)^{*}_{}b(q)}=\gamma_{c(q)
b(q)^{*}}\gamma_{b(q)}=\hbox{\rm Ad }\tau_{c(q)b(q)^{*}}\circ\gamma_{
b(q)}\,,$$
and for $q,r\in S/N$ we have
$$\eqalign{V_q\alpha_q(1_{M(E_{q^{*}})}V_r)u_{q,r}V_{qr}^{*}&=\tau_{
c(q)b(q)^{*}}\gamma_{b(q)}(1_{M(E_{q^{*}})}\tau_{c(r)b(r)^{*}})\tau_{
b(q)b(r)b(qr)^{*}}\tau_{c(qr)b(qr)^{*}}^{*}\cr
&=\tau_{c(q)}1_{M(E_{q^{*}})}\tau_{c(r)b(r)^{*}}\tau_{b(r)b(qr)^{
*}}\tau_{c(qr)b(qr)^{*}}^{*}\cr
&=\tau_{c(q)}1_{M(E_{q^{*}})}\tau_{c(r)}\tau_{b(qr)^{*}}\tau_{b(q
r)^{*}}^{*}\tau_{c(qr)}^{*}=w_{q,r}\,.\cr}
$$
\Eop

\uj Our next goal is to show that starting with a 
Busby-Smith twisted action we can build a 
Green twisted action with the same crossed product.  
The construction mimics that of Fell's in [Fel, Theorem 
I.9.1].  

\lemma Let $(A,T,\beta ,w)$ be a Busby-Smith twisted action. The 
set  
$$S=\{u\partial_t:u\in UM(E_t),\quad t\in T\}$$
with multiplication and adjoint defined by 
$$\eqalign{u_r\partial_r*u_t\partial_t&=\beta_r(\beta_r^{-1}(u_r)
u_t)w_{r,t}\partial_{rt}\cr
(u_t\partial_t)^{*}&=\beta_t^{-1}(u_t^{*})w_{t^{*},t}^{*}\partial_{
t^{*}}\cr}
$$
is an inverse semigroup with idempotent semilattice
$$E_S=\{1_{M(E_f)}\partial_f:f\in E_T\},$$
where $E_T$ is the idempotent semilattice of $T$. If for 
$u_t\partial_t\in S$, $a\in E_{t^{*}}$ and $u_f\partial_f$ in
$$N:=\{u\partial_f:u\in UM(E_f),\quad f\in E_T\}$$
we define 
$$\gamma_{u_t\partial_t}(a)=u_t\beta_t(a)u_t^{*}\quad\hbox{\rm and}
\quad\tau_{u_f\partial_f}=u_f\,,$$
then $(A,S,N,\gamma ,\tau )$ is a Green twisted action.  
\label{connectionlemma2}

\proof To verify that $S$ is an inverse semigroup note that 
the operations are well-defined and as in the proof of 
Proposition \cite{banachref}, associativity holds for the 
multiplication.  Also if $u_t\partial_t\in S$ then by the calculation 
in the proof of Proposition \cite{banachref} and 
Proposition \cite{aproperty}(k) we have 
$$\eqalign{u_t\partial_t(u_t\partial_t)^{*}u_t\partial_t&=\beta_t
(\beta_t^{-1}(u_t)\beta_{t^{*}}(\beta_{t^{*}}^{-1}(\beta_t^{-1}(u_
t^{*})w_{t^{*},t}^{*})u_t))w_{t,t^{*}}w_{tt^{*},t}\partial_{tt^{*}
t}\cr
&=\beta_t(\beta_t^{-1}(u_t)\beta_t^{-1}(u_t^{*})w_{t^{*},t}^{*})\beta_{
t^{*}}(u_t))w_{t,t^{*}}\partial_t\cr
&=u_tu_t^{*}\beta_t(w_{t^{*},t}^{*})\beta_t(\beta_{t^{*}}(u_t))w_{
t,t^{*}}\partial_t\cr
&=w_{t,t^{*}}^{*}\beta_t(\beta_{t^{*}}(u_t))w_{t,t^{*}}\partial_t
=u_t\partial_t\cr}
$$
and similarly $(u_t\partial_t)^{*}u_t\partial_t(u_t\partial_t)^{*}
=(u_t\partial_t)^{*}$. 
It is easy to see that $E_S$ is in fact the idempotent 
semilattice of $S$. The definition
$$u\partial_r\sim v\partial_t\hbox{\rm \ if and only if }r=t$$
gives an idempotent-separating congruence on $S$.  The 
corresponding normal Clifford subsemigroup is $N$.  

If $a\in E_{(rt)^{*}}$ then we have 
$$\eqalign{\gamma_{u_r\partial_r}(\gamma_{u_t\partial_t}(a))&=u_r
\gamma_r(u_t\gamma_t(a)u_t^{*})u_r^{*}=\beta_r(\beta_r^{-1}(u_r)u_
t\beta_t(a)u_t^{*}\beta_r^{-1}(u_r^{*}))\cr
&=\beta_r(\beta_r^{-1}(u_r)u_t)w_{r,t}\beta_{rt}(a)w_{r,t}^{*}\beta_
r(u_t^{*}\beta_r^{-1}(u_r^{*}))\cr
&=\gamma_{\beta_r(\beta_r^{-1}(u_r)u_t)w_{r,t}\partial_{rt}}(a)=\gamma_{
u_r\partial_r*u_t\partial_t}(a)\,,\cr}
$$
showing that $\gamma$ is a homomorphism.
To verify Definition \cite{gtactiondef}(b) note that by 
Proposition \cite{aproperty}(k) for $s=u_t\partial_t$ we have 
$$ss^{*}=\beta_t(\beta_t^{-1}(u_t)\beta_t^{-1}(u_t^{*})w_{t^{*},t}^{
*})w_{t,t^{*}}\partial_{tt^{*}}=\beta_t(w_{t^{*},t}^{*})w_{t,t^{*}}
\partial_{tt^{*}}=1_{M(E_{tt^{*}})}\partial_{tt^{*}},$$
for all $s\in S$.  Thus if $n=u_f\partial_f\in N$ and $s=u_t\partial_
t\in S$ with 
$n^{*}n\leq s^{*}s$ then 
$$1_{M(E_{ft^{*}t})}\partial_{ft^{*}t}=1_{M(E_f)}\partial_f1_{M(E_{
tt^{*}})}\partial_{tt^{*}}=n^{*}ns^{*}s=n^{*}n=1_{M(E_f)}\partial_
f,$$
and so $f\leq t^{*}t$. Hence $\tau_n=u_f\in E_{t^{*}t}=E_t$ and we have
$$\eqalign{\tau_{sns^{*}}&=\beta_t(\beta_t^{-1}(u_t)u_f\beta_t^{-
1}(u_t^{*})w_{t^{*},t}^{*})w_{t,f}w_{tf,t}\cr
&=u_t\beta_t(u_f)u_t^{*}w_{t,t^{*}t}w_{tt^{*},t}^{*}w_{t,t^{*}}^{
*}\cr
&=u_t\beta_t(u_f)u_t^{*}=\gamma_s(\tau_n)\,.\cr}
$$
The conditions in Definition \cite{gtactiondef}(a) and (c) 
are easy to verify. 
\eop

\definition The Busby-Smith twisted actions $(A,S,\alpha ,u)$ and $
(B,T,\beta ,w)$ 
are called {\it conjugate}, if there is a pair $(\rho ,\phi )$ 
such that $\rho :A\to B$ and $\phi :S\to T$ are isomorphisms such that 
for all $s,t\in S$ we have 
$$\rho\circ\alpha_s=\beta_{\phi (s)}\circ\rho\quad\hbox{\rm and}\quad
\rho (u_{s,t})=w_{\phi (s),\phi (t)}\,,$$
where $\rho$ is extended to the double dual of $A$.

\uj
The proof of the following is straightforward.

\lemma Conjugate Busby-Smith twisted actions have isomorphic 
crossed products.
\label{conjugacyref}

\theorem If $(A,T,\beta ,w)$ is a Busby-Smith twisted action and $
(A,S,N,\gamma ,\tau )$ 
is the corresponding Green twisted action constructed 
in Lemma \cite{connectionlemma2}, then the crossed 
products $A\times_{\beta ,w}T$ and $A\times_{\gamma ,\tau}S$ are isomorphic.  
\label{connectionthm2} 

\proof It is easy to see that the cross-section $c:S/N\to S$ 
defined by $c([u\partial_t])=1_{M(E_t)}\partial_t$ is order-preserving so if 
for all $q,r\in S/N$ 
$$\tilde{\beta}_q=\gamma_{c(q)}\quad\hbox{\rm and}\quad\tilde {w}_{
q,r}=\tau_{c(q)*c(r)*c(qr)^{*}}\,,$$
then by Theorem \cite{connectionthm1}, $(A,S/N,\tilde{\beta },\tilde {
w})$ 
is a Busby-Smith twisted action and the crossed 
products $A\times_{\gamma ,\tau}S$ and $A\times_{\tilde{\beta },\tilde {
w}}S/N$ are isomorphic.  We 
finish the proof by showing that $(\hbox{\rm id},\phi )$ is a conjugacy 
between the Busby-Smith twisted actions $(A,T,\beta ,w)$ and 
$(A,S/N,\tilde{\beta },\tilde {w})$, where 
$$\phi (t)=[1_{M(E_t)}\partial_t].$$
It is easy to see that $\phi$ is an isomorphism. For $s,t\in T$ 
we have 
$$\tilde{\beta}_{\phi (s)}=\gamma_{c([1_{M(E_s)}\partial_s])}=\beta_
s$$
and
$$\eqalign{\tilde {w}_{\phi (s),\phi (t)}&=\tau_{c([1_{M(E_s)}\partial_
s])c([1_{M(E_t)}\partial_t])c([1_{M(E_{st})}\partial_{st}])^{*}}\,\cr
&=\tau_{1_{M(E_s)}\partial_s*1_{M(E_t)}\partial_t*(1_{M(E_{st})}\partial_{
st})^{*}}\,\cr
&=\beta_{qr}(\beta_{qr}^{-1}(w_{q,r})w_{(qr)^{*},qr}^{*})w_{qr,(q
r)^{*}}=w_{q,r}\,,\cr}
$$
and so we are done by Lemma \cite{conjugacyref}.
\eop

\newsection{Decomposition of Green twisted actions}

\uj Now we prove that in the presence of a normal 
Clifford subsemigroup, a Green twisted crossed product 
can be decomposed as an iterated Green twisted crossed 
product.  The close analogy with Green's decomposition 
theorem [Gre] gives further evidence that in the inverse 
semigroup case, normal Clifford subsemigroups play the 
same role as normal subgroups do in the group case.  

\theorem Let $(A,S,N,\gamma ,\tau )$ be a Green twisted  
action and $K$ be a normal Clifford subsemigroup 
containing $N$.  The restriction of $\gamma$ to $K$ gives a Green 
twisted action $(A,K,N,\gamma ,\tau )$.  Let $\mu\times_{\tau}m$ be the 
universal representation of $A\times_{\gamma ,\tau}K$, and identify $
A\times_{\gamma ,\tau}K$ 
with its image under $\mu\times_{\tau}m$.  Let 
$$\widetilde{E}_s=\overline {\hbox{\rm span}}\,\{\mu (a)m_k:a\in E_
k,k\in K,kk^{*}\leq ss^{*}\}$$
and define 
$$\tilde{\gamma}_s:\widetilde{E}_{s^{*}}\to\widetilde{E}_s\quad\hbox{\rm by}
\quad\tilde{\gamma}_s=\hbox{\rm Ad}\,m_s\,.$$
Let $\tilde{\tau}_k=m_k$ for all $k\in K$.  Then $(A\times_{\gamma 
,\tau}K,S,K,\tilde{\gamma },\tilde{\tau })$ is a 
Green twisted action and 
$$A\times_{\gamma ,\tau}S\cong (A\times_{\gamma ,\tau}K)\times_{\tilde{
\gamma },\tilde{\tau}}S\,.$$

\label{greendecomp}

\proof First we show that $(\tilde{\gamma },\tilde{\tau })$ is a Green twisted 
action.  $\widetilde E_s$ is an ideal by the calculation in the 
proof of Proposition \cite{hivcpv} and the fact that 
$k,l\in K$ and $kk^{*}\leq ss^{*}$ imply $kl(kl)^{*}\leq ss^{*}$ and $
lk(lk)^{*}\leq ss^{*}$.  It is 
clear that $\tilde{\gamma}_s$ is multiplicative and preserves adjoints.  
Since 
$$\tilde{\gamma}_s(\mu (a)m_k)=m_s\mu (a)m_km_s^{*}=m_s\mu (a)m_s^{
*}m_{sks^{*}}=\mu (\gamma_s(a))m_{sks^{*}}$$
for all $\mu (a)m_k\in\widetilde E_k$, $\tilde{\gamma}_s$ is a bijection and 
therefore an isomorphism of $\widetilde E_{s^{*}}$ to 
$\widetilde E_s$, hence a partial automorphism of $A\times_{\gamma 
,\tau}K$.  For 
$\mu (a)m_k\in\widetilde E_{(st)^{*}}$ we have 
$$\tilde{\gamma}_s\tilde{\gamma}_t(\mu (a)m_k)=\tilde{\gamma}_s(\mu 
(\gamma_t(a))m_{tkt^{*}})=\mu (\gamma_s(\gamma_t(a))m_{stkt^{*}s^{
*}}=\tilde{\gamma}_{st}(\mu (a)m_k)\,,$$
so to show that $\tilde{\gamma}$ is a homomorphism from $S$ to 
$\hbox{\rm PAut}\,A\times_{\gamma ,\tau}K$ we only need to check that 
$\hbox{\rm dom}\,\tilde{\gamma}_s\tilde{\gamma}_t=\hbox{\rm dom}\,
\tilde{\gamma}_{st}$.  First note that if $T\in\widetilde A_t$ 
then $T=\lim_iT_i$ for some $T_i\in\hbox{\rm span}\,\{\mu (a)m_k:
a\in A_k,k\in K,kk^{*}\leq tt^{*}\}$.  
Hence, $\mu (1_{M(A_t)})T=\lim_i\mu (1_{M(A_t)})T_i=T$ since for $
a\in A_k$, $k\in K$ and 
$kk^{*}\leq tt^{*}$ we have $A_k\subset A_t$ and so 
$$\mu (1_{M(A_t)})\mu (a)m_k=\mu (1_{M(A_t)}a)m_k=\mu (a)m_k\,.$$
Now, if $T\in\hbox{\rm dom}\,\tilde{\gamma}_s\tilde{\gamma}_t$ then $
\tilde{\gamma}_t(T)=\lim_iT_i$ for some 
$T_i\in\hbox{\rm span}\,\{\mu (a)m_k:a\in A_k,k\in K,kk^{*}\leq s^{
*}s\}$. Then 
$$T=m_t^{*}\mu (1_{M(A_t)})\tilde{\gamma}_t(T)m_t=\lim_im_t^{*}\mu 
(1_{M(A_t)})T_im_t\,,$$
so to see that $T\in\hbox{\rm dom}\,\tilde{\gamma}_{st}$, it suffices to show that if $
a\in A_k$, $k\in K$ and 
$kk^{*}\leq s^{*}s$, then 
$$\eqalign{m_t^{*}\mu (1_{M(A_t))}\mu (a)m_km_t=\mu (\gamma_{t^{*}}
(1_{M(A_t)}a))m_{t^{*}kt}\in\widetilde A_{(st)^{*}}\,.\cr}
$$
This is true since $\gamma_{t^{*}}(1_{M(A_t)}a)\in A_{t^{*}k}=A_{
t^{*}kk^{*}t}$ and 
$t^{*}kk^{*}t\leq t^{*}s^{*}st=(st)^{*}st$.  

To see that $\hbox{\rm dom}\,\tilde{\gamma}_s\tilde{\gamma}_t\supset\hbox{\rm dom}\,
\tilde{\gamma}_{st}$, first note that 
$\widetilde A_{t^{*}}\supset\widetilde A_{(st)^{*}}$ since $kk^{*}
\leq (st)^{*}st$ implies 
$kk^{*}\leq t^{*}s^{*}st\leq t^{*}t$. So it suffices to show that if $
a\in A_k$, $k\in K$ and 
$kk^{*}\leq (st)^{*}st$ then 
$$\tilde{\gamma}_t(\mu (a)m_k=m_t\mu (a)m_km_t^{*}=\mu (\gamma_t(
a))m_{tkt^{*}}\in\widetilde A_{s^{*}}\,.$$
But this follows since 
$$\gamma_t(a)\in A_{tk}=A_{tkk^{*}t^{*}}=A_{tkt^{*}tk^{*}t^{*}}=A_{
tkt^{*}(tkt^{*})^{*}}=A_{tkt^{*}}$$
and $tkt^{*}(tkt^{*})^{*}=tkk^{*}t^{*}\leq tt^{*}s^{*}stt^{*}\leq 
s^{*}s$.  We used the fact that 
$tkk^{*}t^{*}=tkt^{*}tk^{*}t^{*}$.  This is true because $tkk^{*}
t^{*}$ is the unique 
idempotent in $[t][kk^{*}][t^{*}]$ and $tkt^{*}tk^{*}t^{*}$ is the unique idempotent in 
$[t][k][t^{*}t][k^{*}][t^{*}]=[t][kk^{*}][t^{*}t][kk^{*}][t^{*}]=
[t][kk^{*}][t^{*}]$.  

Next we show that $\tilde{\tau}_k$ is a unitary multiplier of 
$\widetilde E_k$.  If $a\in E_l$ for some $l\in K$ with $ll^{*}\leq 
kk^{*}$ 
then we have $\mu (a)m_lm_k=\mu (a)m_{lk}\in\widetilde E_k$ since 
$lk(lk)^{*}=ll^{*}$ and so $lk(lk)^{*}\leq kk^{*}$ and 
$a\in E_l=E_{ll^{*}}=E_{lk(lk)^{*}}=E_{lk}$.  We also have 
$m_k\mu (a)m_l=(\mu (\tau_l^{*}a^{*})m_{k^{*}})^{*}\in\widetilde 
E_k$ since $\tau_l^{*}a^{*}\in E_l$ and 
so $\tau_l^{*}a^{*}\in E_l=E_{kk^{*}l}=\gamma_{kk^{*}}(E_{kk^{*}}
E_l)\subset E_{kk^{*}}=E_k$.  Hence 
$\tilde{\tau}_k$ is a multiplier of $\widetilde E_k$.  The multiplier $
m_k$ is 
clearly unitary with inverse $m_k^{*}$.  

To check Definition \cite{gtactiondef}(a) let $k\in K$ and 
$\mu (a)m_l\in\widetilde E_{k^{*}}$. Then $a\in E_l$ for some $l\in 
K$ such that 
$l^{*}$$l\leq k^{*}k$, so 
$$\tilde{\gamma}_k(\mu (a)m_l)=m_k\mu (a)m_lm_k^{*}=\tilde{\tau}_
k\mu (a)m_l\tilde{\tau}_k^{*}\,.$$
Since for $s\in S$ and $k\in K$ with $k^{*}k\leq s^{*}s$ we have 
$$\tilde{\gamma}_s(\tilde{\tau}_k)=\tilde{\gamma}_s(m_k)=m_sm_km_
s^{*}=m_{sks^{*}}=\tilde{\tau}_{sks^{*}}\,,$$
Definition \cite{gtactiondef}(b) is verified and so we are 
finished showing that $(\tilde{\gamma },\tilde{\tau })$ is a Green twisted action.

Next we find a correspondence between covariant 
representations of our actions.
First let $(\Pi ,v)$ be a covariant representation of $(\tilde{\gamma }
,\tilde{\tau }).$ 
Then $\Pi =\pi\times_{\tau}u$ for some covariant representation $
(\pi ,u)$ of 
$(\gamma |_K,\tau )$. We show that $(\pi ,v)$ is a covariant representation 
of $(\gamma ,\tau )$ by checking the conditions of Definition 
\cite{gcovrepdef}. 
Condition (a) follows from the calculation
$$\eqalign{v_s\pi (a)v_s^{*}&=v_s\Pi (\mu (a)m_{s^{*}s})v_s^{*}=\Pi 
(\tilde{\gamma}_s(\mu (a)m_{s^{*}s}))\cr
&=\Pi (\gamma_s(a)m_{ss^{*}ss^{*}})=\pi (\gamma_s(a))\,.\cr}
$$
Since $v_s$ has final space 
$$\eqalign{\Pi (\widetilde{E}_s)H&=\overline {\hbox{\rm span}}\,\{
\Pi (\mu (E_k)m_k):k\in K,kk^{*}\leq ss^{*}\}\cr
&=\Pi (\mu (E_{ss^{*}})m_{ss^{*}})=\pi (E_s)\,\,,\cr}
$$
$v_s$ has the required initial and final spaces. The 
equality $v_n=\Pi (\tilde{\tau}_n)=\Pi (m_n)=\Pi (\mu (\tau_n))=\pi 
(\tau_n)$ for all 
$n\in N$ shows that 
$(\pi ,v)$ preserves the twist.  

Next we show that if $(\pi ,v)$ is a covariant representation of $
(\gamma ,\tau )$,
then $(\pi\times_{\tau}v|K,v)$ is a covariant representation of 
$(\tilde{\gamma },\tilde{\tau })$. Condition (a) follows from the calculation
$$\eqalign{v_s\big(\pi\times_{\tau}v|K\big)(\mu (a)m_k)v_s^{*}&=v_
s\pi (a)v_kv_s^{*}=v_s\pi (a)v_s^{*}v_sv_kv_s^{*}\cr
&=\pi (\gamma_s(a))v_{sks^{*}}=\big(\pi\times_{\tau}v|K\big)(\mu 
(\gamma_s(a))m_{sks^{*}})\cr
&=\big(\pi\times_{\tau}v|K\big)(\tilde{\gamma}_s(\mu (a)m_k))\,\,
,\cr}
$$
where $a\in E_k$ for some $k\in K$ satisfying $kk^{*}\leq ss^{*}$.
Conditions (b) and (c) are clearly satisfied. 

Now we investigate the ranges of the corresponding 
representations. If $(\pi\times_{\tau}u,v)$ is a covariant representation 
of $(\tilde{\gamma },\tilde{\tau })$ and $(\pi ,v)$ is a covariant representation of 
$(\gamma ,\tau )$ then for $a\in E_s$ and $s\in S$ we have
$$\big(\pi\times_{\tau}v\big)(\mu (a)m_s)=\pi (a)v_s=\pi (a)u_{ss^{
*}}v_s=\big((\pi\times_{\tau}u)\times_{\tilde{\tau}}v\big)(\mu (a)m_{
ss^{*}}\delta_s)\,.$$
Hence the range of $\pi\times_{\tau}v$ is contained in the range of 
$(\pi\times_{\tau}u)\times v$ since $\mu (a)m_{ss^{*}}\in\widetilde 
E_s$. On the other hand, 
for $\mu (a)m_k\in\widetilde E_{s^{*}}$, that is, $a\in E_k$ for some $
k\in K$ with 
$kk^{*}\leq ss^{*}$, we have  
$$\big((\pi\times_{\tau}u)\times_{\tilde{\tau}}v\big)(\mu (a)m_k\delta_
s)=\pi (a)\big(\pi\times_{\tau}u\big)(\tilde{\tau}_k)v_s=\pi (a)v_
kv_s=\big(\pi\times_{\tau}v\big)(\mu (a)m_{ks})\,.$$
Hence the range of $(\pi\times_{\tau}u)\times v$ is contained in the range of 
$\pi\times_{\tau}v$.

Now starting with a faithful representation $\rho\times_{\tau}z$ of 
$A\times_{\gamma ,\tau}S$ we know that $(\rho\times_{\tau}z|K)\times_{
\tilde{\tau}}z$ is a representation of 
$(A\times_{\gamma ,\tau}K)\times_{\tilde{\gamma },\tilde{\tau}}S$ and $
\Psi =(\rho\times_{\tau}z)^{-1}\circ (\rho\times_{\tau}z|K)\times_{
\tilde{\tau}}z$ is a 
surjective homomorphism.  Similarly, starting with a 
faithful representation $(\pi\times_{\tau}u)\times_{\tilde{\tau}}v$ of $
(A\times_{\gamma ,\tau}K)\times_{\tilde{\gamma },\tilde{\tau}}S$ we 
know that $\pi\times_{\tau}v$ is a representation of $A\times_{\gamma 
,\tau}S$ and 
$\Phi =((\pi\times_{\tau}u)\times_{\tilde{\tau}}v)^{-1}\circ (\pi\times_{
\tau}v)$ is a surjective homomorphism.  
The relationship among these maps is expressed by the 
commutative diagram:
 
$$\matrix{A\times_{\gamma ,\tau}S&\buildrel\rho\times_\tau z \over\mapeitoionto&\hbox{\rm im}\,(\rho\times_\tau z)\cr
\cr
\scriptstyle \pi\times_\tau v \mapsonto\quad&& \qquad\mapnonto\scriptstyle (\rho\times_\tau z|K)\times_{\tilde\tau} z\cr
\cr
\hbox{\rm im}\,((\pi\times_\tau u)\times_{\tilde\tau} v)&\mathop{\mapwitoionto}\limits_{(\pi\times_\tau u)\times_{\tilde\tau} v} &(A\times_{\gamma ,\tau}K)\times_{\tilde{\gamma }
,\tilde{\tau}}S\cr
\cr}
$$
We show that $\Psi\circ\Phi$ is the identity map, which 
implies that $\Phi$ is an isomorphism. For $a\in E_s$ we 
have 
$$\eqalign{\Psi\circ\Phi (a\delta_s)&=\Psi\circ ((\pi\times_{\tau}
u)\times_{\tilde{\tau}}v)^{-1}(\pi (a)v_s)=\Psi (\mu (a)m_{ss^{*}}\delta_
s)\cr
&=(\rho\times_{\tau}z)^{-1}(\rho (a)z_{ss^{*}}z_s)=\mu (a)m_s\,.\cr}
\,$$
\Eop

\uj Recall [Sie] that if $\beta$ is the canonical action of the 
inverse semigroup $S$ on its semilattice $E$ (Example 
\cite{canonicalaction}), then $C^{*}(S)\cong C^{*}(E)\times_{\beta}
S$.  Also 
recall that if $\beta$ is the trivial action of $S$ on {\bf C}, that is, 
$\beta_s$ is the identity for all $s\in S$, then $C^{*}(G_S)\cong 
{\bf C}\times_{\beta}S$, 
where $G_S$ is the maximal group homomorphic image of 
$S$.  These results can be extended to Green twisted 
inverse semigroup actions.  Similar results will be given 
for Busby-Smith twisted actions in Proposition 
\cite{tactionexproof}.  

\prop Let $S$ be an inverse semigroup  
with idempotent semilattice $E$, and let $K$ be a normal 
Clifford subsemigroup of $S$.  
For 
all $s\in S$ let $E_s$ be the closed span of $\{k\in K:kk^{*}\leq 
ss^{*}\}$ in 
$C^{*}(K)$, and define $\alpha_s:E_{s^{*}}\to E_s$ 
by $\alpha_s(a)=sas^{*}$. For $k\in K$ let $\sigma_k=k$.  
Then $(C^{*}(K),S,K,\alpha ,\sigma )$ is a Green twisted action 
and $C^{*}(S)\cong C^{*}(K)\times_{\alpha ,\sigma}S$.
\label{btactionproof}

\proof Let $\beta$ be the canonical action of $S$ on its 
semilattice $E$.  The result can be obtained by following 
through the isomorphism chain 
$$C^{*}(S)\cong C^{*}(E)\times_{\beta}S\cong (C^{*}(E)\times_{\beta}
K)\times_{\tilde{\beta },\tilde{\tau}}S\cong C^{*}(K)\times_{\alpha 
,\sigma}S,$$
using Theorem \cite{greendecomp}.  \Eop

\prop Let $S$ be an inverse semigroup with idempotent 
semilattice $E$, and let $K$ be a normal Clifford 
subsemigroup of $S$.  For all $s\in S$ let $E_s$ be the closed 
span in $C^{*}(G_K)$ of $\{[k]:k\in K,kk^{*}\leq ss^{*}\}$, where $
G_K$ 
is the maximal group homomorphic image of $K$ and $[k]$ is 
the canonical image of $k$ in $C^{*}(G_K)$, and define 
$\alpha_s:E_{s^{*}}\to E_s$ by $\alpha_s(a)=sas^{*}$.  For $k\in 
K$ let $\sigma_k=[k]$.  
Then $(C^{*}(G_K),S,K,\alpha ,\sigma )$ is a Green twisted action and 
$C^{*}(G_S)\cong C^{*}(G_K)\times_{\alpha ,\sigma}S$, where $G_S$ is the maximal group 
homomorphic image of $S$.  
\label{btactionex2} 
 
\proof Let $\beta$ be the trivial action of $S$ on {\bf C}.  The result 
can be obtained by following through the isomorphism 
chain 
$$C^{*}(G_S)\cong {\bf C}\times_{\beta}S\cong ({\bf C}\times_{\beta}
K)\times_{\tilde{\beta },\tilde{\tau}}S\cong C^{*}(K)\times_{\alpha 
,\sigma}S,$$
using Theorem \cite{greendecomp}.  \Eop

\uj Recall from [Re1] that the Cuntz algebra is the $C^{*}$-algebra $
C^{*}(O_n)$ 
of the Cuntz groupoid $O_n$.  Paterson recently found a connection 
between groupoid $C^{*}$-algebras and inverse semigroup crossed products 
[Pa3].  The details of this connection will be included in his upcoming 
book.  Using his results it is possible to identify $C^{*}(O_n)$ with a 
crossed product of $C_0(O_n^0)$ by the Cuntz inverse semigroup ${\cal O}_
n$.  

Since the only normal Clifford subsemigroup of the Cuntz inverse 
semigroup is its idempotent semilattice, there is no nontrivial 
decomposition of this crossed product.  This is further evidence of 
the rigidity of the Cuntz relations.  The question remains whether 
the Cuntz-Krieger inverse semigroups [HR] have normal Clifford 
subsemigroups supplying a possible decomposition of the 
Cuntz-Krieger algebras.

\newsection{Decomposition of Busby-Smith twisted actions} 

\uj Since every Busby-Smith twisted action corresponds to a Green 
twisted action, we can apply our Decomposition Theorem 
\cite{greendecomp}
to Busby-Smith twisted actions.  

\theorem If $(A,T,\beta ,w)$ is a Busby-Smith twisted action and $
L$ 
is a normal Clifford subsemigroup of $T$ with an 
order-preserving cross-section $c:T/L\to T$, then the 
crossed product can be decomposed as an iterated 
Busby-Smith twisted crossed product 
$$A\times_{\beta ,w}T\cong (A\times_{\beta ,w}L)\times_{\tilde{\beta }
,\tilde {w}}T/L.$$

\label{busbydecomp}

\proof Let $(A,S,N,\gamma ,\tau )$ be the Green twisted action 
corresponding to $(A,T,\beta ,w)$, constructed in Lemma 
\cite{connectionlemma2}.  Let $\sim_L$ denote the 
idempotent-separating congruence on $T$ determined by $L$.  
Define a normal Clifford subsemigroup 
$$K=\{u_t\partial_t:u_t\in UM(E_t),\quad t\in L\}$$
of $S$. The corresponding congruence relation on $S$ is 
given by 
$$u\partial_r\sim_Kv\partial_t\hbox{\rm \ if and only if }r\sim_L
t.$$
By Theorem \cite{greendecomp} there is a Green twisted 
action $(A\times_{\gamma ,\tau}K,S,K,\tilde{\gamma },\tilde{\tau }
)$ such that $A\times_{\gamma ,\tau}S$ 
is isomorphic to $(A\times_{\gamma ,\tau}K)\times_{\tilde{\gamma }
,\tilde{\tau}}S$. Define a cross-section 
$$d:S/K\to S,\qquad d([u_t\partial_t])=1_{M(E_{c([t])})}\partial_{
c([t])}\,.$$
To see that $d$ is order-preserving, let
$[u_s\partial_s]\leq [u_t\partial_t]$ in $S/K$. Then $[u_s\partial_
s]=[u_t\partial_t*(u_s\partial_s)^{*}*u_s\partial_s]$ 
and so $s\sim_Lts^{*}s$. This means $[s]\leq [t]$ in $T/L$ and so 
$d([u_s\partial_s])\leq d([u_t\partial_t])$ because 
$$\eqalign{d([u_t\partial_t])d([u_t\partial_t])^{*}d([u_t\partial_
t])&=1_{M(E_{c([t])})}\partial_{c([t])}*1_{M(E_{c([t])^{*}})}\partial_{
c([t])^{*}c([t])}\cr
&=w_{c([t]),c([s])^{*}c([s])}\partial_{c([t])c([s])^{*}c([s])}\cr
&=1_{M(E_{c([s])})}\partial_{c([s])}\,.\cr}
$$
Hence by Theorem 
\cite{connectionthm1} there is a Busby-Smith twisted  
action $(A\times_{\gamma ,\tau}K,S/K,\tilde{\alpha },\tilde {u})$ such that $
(A\times_{\gamma ,\tau}K)\times_{\tilde{\gamma },\tilde{\tau}}S$ and 
$(A\times_{\gamma ,\tau}K)\times_{\tilde{\alpha },\tilde {u}}S/K$ are isomorphic.  By Theorem 
\cite{connectionthm2} there is an isomorphism 
$\rho :A\times_{\beta ,w}L\to A\times_{\gamma ,\tau}K$, and it is easy to see that 
$$\phi :T/L\to S/K,\qquad\phi ([t])=[1_{M(E_t)}\partial_t]$$
is an isomorphism.  So if for all $s,t\in T/L$ we define 
 
$$\tilde{\beta}_s=\rho^{-1}\circ\tilde{\alpha}_{\phi (s)}\circ\rho
\quad\hbox{\rm and}\quad\tilde {w}_{s,t}=\rho^{-1}(\tilde {u}_{\phi 
(s),\phi (t)})\,,$$
then $(A\times_{\beta ,w}L,T/L,\tilde{\beta },w)$ is clearly a Busby-Smith twisted 
action, 
which is conjugate to $(A\times_{\gamma ,\tau}K,S/K,\tilde{\alpha }
,\tilde {u})$.  Thus we have 
the following chain of isomorphisms 
$$\eqalign{A\times_{\beta ,w}T&\cong A\times_{\gamma ,\tau}S\cr
&\cong (A\times_{\gamma ,\tau}K)\times_{\tilde{\gamma },\tilde{\tau}}
S\cr
&\cong (A\times_{\gamma ,\tau}K)\times_{\tilde{\alpha },\tilde {u}}
S/K\cr
&\cong (A\times_{\beta ,w}L)\times_{\tilde{\beta },\tilde {w}}T/L\,
.\cr}
$$
\Eop

\uj Using Theorem \cite{connectionthm1} we get analogs of 
Propositions \cite{btactionproof} and \cite{btactionex2} for 
Busby-Smith twisted actions:
 
\prop If $(C^{*}(N),T/N,\beta ,w)$ is the Busby-Smith twisted 
action of Proposition \cite{tactionexample}, then
$$C^{*}(T)\cong C^{*}(N)\times_{\beta ,w}T/N\,.$$
\label{tactionexproof} 

\prop Let $S$ be an inverse semigroup and let $N$ be a normal 
Clifford subsemigroup with an order-preserving 
cross-section from $S/N$ to $S$.  Then
$$C^{*}(G_S)\cong C^{*}(G_N)\times_{\beta ,w}S/N\,.$$

\uj
\noindent This material is based upon work supported
by the National Science Foundation under Grant No.
DMS9401253.

\bigskip

\centerline{\nagy References}
\nobreak
\medskip
\nobreak

\def\i#1#2{\itemitem {\hbox to .4in{#1\hfil}} {#2} }
\rm

\i{[BS]}{R.C. Busby and H.A. Smith, \it Representation of twisted group
algebras\rm, Trans. Amer. Math. Soc. \bf 149 \rm (1970), 503--537.}
\i{[CP]}{A. Clifford and G. Preston, {\it The algebraic theory of 
semigroups}, American Mathematical Society, Providence, 
Rhode Island, 1961.}
\i{[DP1]}{J. Duncan, S.L.T. Paterson, $C^{*}${\it -algebras of inverse semigroups},
Proc. Edinburgh Math. Soc.
{\bf 28} (1985), 41-58.}
\i{[DP2]}{J. Duncan, S.L.T. Paterson, $C^{*}${\it -algebras of Clifford }
{\it semigroups}, Proc. Roy{\it .\/} Soc. Edinburgh, {\bf 111}(A) (1989),
129-145.}
\i{[Ex1]} {R. Exel, {\it Circle actions on $C^{*}$-algebras, 
partial automorphisms
and a generalized Pimsner-Voiculescu exact sequence}, J.  
Funct.  Anal. {\bf 122}  (1994), 361-401.}
\i{[Ex2]}{R. Exel, {\it Twisted partial actions a classification 
of stable $C^{*}$-algebraic bundles}, preprint.}
\i{[Ex3]}{R. Exel, \it Partial actions of groups and actions 
of inverse semigroups\rm, preprint.}
\i{[Fel]}{J.M.G. Fell, {\it An extension of Mackey's method to 
Banach $*$-algebraic bundles}, Memoirs Amer. Math. Soc. 90, Providence,
Rhode Island, 1969.}
\i{[Gre]} {P. Green, {\it The local structure of twisted 
covariance algebras}, Acta Math. \bf 140 \rm (1978), 191-250.}
\i{[How]} {J.M. Howie, {\it An introduction to semigroup theory}, 
Academic Press, London, 1976.}
\i{[Kum]}{A. Kumjian, {\it On localization and simple }
$C^{*}${\it -algebras}, Pacific J. Math. \bf 112 \rm(1984), 141-192.}
\i{[HR]}{R. Hancock, I. Raeburn, {\it The $C^*$-algebras of some inverse
semigroups}, Bull. Austral. Math. Soc. {\bf 42} (1990), 335-348.}
\i{[McC]} {K. McClanahan, \it K-theory for partial crossed products by discrete
 groups\rm, J. Funct. Anal. {\bf 130} (1995), 77-117.} 
\i{[Nic]}{A. Nica, {\it On a groupoid construction for actions
of certain inverse semigroups}, preprint.}
\i{[PR]}{J. Packer and I. Raeburn, {\it Twisted crossed }
{\it products of} $C^{*}${\it -algebras}, Math. Proc. Camb. Phil. Soc. 
\bf 106 \rm(1989), 293-311.} 
\i{[Pa1]}{A.L.T. Paterson, {\it Weak containment and Clifford }
{\it semigroups}, Proc. Roy. Soc. Edinburgh {\bf 81}(A) (1978), 23-30.}
\i{[Pa2]}{A.L.T. Paterson, {\it Inverse semigroups, groupoids and a }
{\it problem of J. Renault}, Algebraic methods in Operator
Theory, ed. R. Curto and P.E.T. Jorgensen, Birkh\"auser,
Boston, 1993, 11 pp.}
\i{[Pa3]}{A.L.T. Paterson, {\it r-Discrete $C^{*}$-algebras as covariant 
$C^{*}$-algebras}, Groupoid Fest lecture notes, Reno, 1996.}
\i{[Pet]}{M. Petrich, {\it Inverse semigroups}, John Wiley \& 
Sons, New York, 1984.}
\i{[Re1]}{J.N. Renault, {\it A groupoid approach to} $C^{*}${\it -algebras},
Lecture Notes in Mathematics, Vol 793, Springer-Verlag,
New York, 1980.}
\i{[Re2]}{J.N. Renault, {\it Repr\'esentation des produits crois\'es }
{\it d'alg\`ebres de groupoides}, J. Operator Theory, {\bf 18} (1987),
67-97.}
\i{[Sie]}{N. Sieben, $C^{*}${\it-crossed products by partial 
actions and actions of inverse semigroups}, J. Australian 
Math. Soc. (Series A) {\bf 63} (1997), 32-46.}

}
\bye